\documentclass{article}
\usepackage[colorlinks]{hyperref}
\usepackage{amsmath,amssymb}
\usepackage{physics}
\usepackage{subcaption}
\usepackage{mathtools}
\usepackage{parskip}
\usepackage{graphicx}
\usepackage{comment}
\usepackage{bm}
\usepackage{amsthm}
\usepackage{siunitx}
\usepackage{upgreek}
\usepackage{tikz}
\usetikzlibrary{calc}
\usepackage[toc, page, title]{appendix}
\usepackage{algorithm}
\usepackage{algpseudocode}
\usepackage{authblk}
\usepackage{braket}
\usepackage[utf8]{inputenc}
\usepackage[pagewise]{lineno}
\usepackage[
	backend=biber,
	style=numeric,
	sorting=none
]{biblatex}

\usetikzlibrary{patterns,arrows,decorations.pathreplacing,math,decorations.markings,angles}
\usepackage[margin=1.0in]{geometry}
\providecommand{\keywords}[1]
{
	\small	
	\textbf{Keywords. } #1
}

\providecommand{\classification}[1]
{
	\small	
	\textbf{MSC subject classifications. } #1
}

\newcommand{\maxw}{\mathcal{M}(f)}
\DeclareMathOperator*{\argmin}{arg\,min}
\DeclarePairedDelimiter\ceil{\lceil}{\rceil}

\newcommand{\transp}{{\tt transp}}

\newtheorem{remark}{Remark}

\newcommand{\subsubsubsection}[1]{\paragraph{#1}\mbox{}\\}
\title{A meshless MUSCL method for the BGK-Boltzmann equation}

\author[1, 2]{Klaas Willems\thanks{Corresponding author: klaas.willems@kuleuven.be}}
\author[2]{Axel Klar\thanks{axel.klar@rptu.de}}
\author[3]{Giovanni Russo\thanks{giovanni.russo1@unict.it}}
\author[1]{Giovanni Samaey\thanks{giovanni.samaey@kuleuven.be}}
\author[2]{Sudarshan Tiwari\thanks{tiwari@mathematik@uni-kl.de}}
\affil[1]{Department of Computer Science, KU Leuven, Leuven, Belgium}
\affil[2]{Faculty of Mathematics, RPTU Kaiserslautern-Landau, Kaiserslautern, Germany}
\affil[3]{Department of Mathematics and Computer Science, University of Catania, Catania 95125, Italy}

\begin{document}
	
\maketitle

\begin{abstract}
We present a numerical method for simulating rarefied gases that interact with moving boundaries and rigid bodies. The gas is described by the BGK equation in Lagrangian form and solved using an Arbitrary Lagrangian-Eulerian method, in which grid points move with the local mean velocity of the gas \cite{tiwari_meshfree_2022}. The main advantage of the moving grid is that the algorithm can deal well with cases where the domain boundaries are time-dependent and the simulation domain contains rigid objects. Due to the irregular nature of the grid, we use a meshless MUSCL-like Moving Least Squares Method (MLS) for spatial discretisation coupled with a higher-order Implicit-Explicit Runge-Kutta method. To avoid spurious oscillations at discontinuities, we use the so-called Multi-dimensional Optimal Order Detection (MOOD) method \cite{clain_high-order_2011} with an adapted criterion to relax the discrete maximum property. Finally, we employ a new implementation of the boundary conditions that requires no iterative or extrapolation procedure. The method achieves fourth-order in 1D and second-order in 2D for simulations with moving boundaries. We demonstrate the method's effectiveness on classical test cases such as the driven square cavity, shear layer, and shock tube.
\end{abstract}
 \\
\classification{74F10, 65M99, 70E99, 35L02}

\section{Introduction}
\label{section:introduction}
The field of rarefied gas dynamics concerns the study of gases in regimes where the mean free path of gas particles is comparable to a characteristic macroscopic length scale of the system. In this regime, classical fluid dynamical models break down and a kinetic description is required. In the case of rarefied gases, the appropriate kinetic model is the Boltzmann equation, a seven-dimensional kinetic integro-differential equation \cite{cercignani_boltzmann_1988}. Historically, the main application of rarefied gas dynamics was in aerospace engineering, where flow is characterized by high Mach numbers and the mean free path of gas molecules in the outer atmosphere is equal to or larger than the size of man-made objects. Simulations in this field are usually performed with the Direct Simulation Monte Carlo method (DSMC) \cite{bird_molecular_1994}, a stochastic particle method for the Boltzmann equation. More recently, the field of rarefied gas dynamics has found applications in micro-electromechanical systems (MEMS) \cite{karniadakis_microflows_2005}. In these applications, a gas interacts with small sensors such as accelerometers, micro pumps and micro engines whose characteristic size is comparable to the mean free path of the gas. The simulation of MEMS is a truly multi-physics problem that requires modelling of structural dynamics, electrical properties and gas dynamics \cite{russo_semilagrangian_2009}. Especially relevant for the gas dynamics is that the gas flows in MEMS require precise treatment of solid-boundary motion. Furthermore, they are characterized by low Mach numbers, therefore DSMC simulations require a significant amount of particles to be simulated to reduce statistical noise. In \cite{zabelok_adaptive_2015}, a parallelized DSMC method for moving boundaries was combined with a fluid solver to deal with this issue. Alternatively, deterministic, grid-based, simulation methods for the Boltzmann equation can be used, of which the most well-known is the spectral method from \cite{filbet_solving_2006, pareschi_fourier_1996, pareschi_numerical_2000}.
On the other hand, a significant amount of deterministic numerical methods rely on a simplified model known as the BGK-Boltzmann equation \cite{bhatnagar_model_1954}. In this model, the collision integral of the Boltzmann equation is replaced by a relaxation towards the steady-state solution. This model is cheaper to work with and considered good enough for flows that are not too rarefied. Numerical methods based on the BGK equation were first applied to MEMS in \cite{frangi_application_2007, frangi_analysis_2008}.

Several classes of numerical methods for the BGK-Boltzmann equation to treat low-Mach flows and moving boundaries have been published. A first class of methods are the so-called semi-Lagrangian methods \cite{santagati_new_2011}. These schemes use the method of characteristics to treat the transport term in the BGK equation. Since the transport is linear, its treatment reduces to an interpolation problem of the distribution function at the previous time step. The BGK relaxation term is treated implicitly to avoid a time step restriction due to the stiffness. However, due to the special structure of the BGK operator, the implicit treatment can be explicitly solved, leading to a very stable and efficient treatment of the collision term. These schemes have been extended to higher-orders \cite{santagati_new_2011, groppi_high_2014}, strict conservation when using nonlinear reconstructions \cite{cho_conservative_2021, cho_conservative_2021-1}, and do not have a time step restriction. Semi-Lagrangian methods have been extended to moving boundary problems in 1D with one-way coupling from the motion of the rigid body to the gas \cite{russo_semilagrangian_2009}. The method relies on a so-called immersed boundary method: the Cartesian grid does not change during computation and the equation is only solved in the grid points in the gaseous domain. Since the exact location of the gas-rigid body is not tracked, this method seems unpractical to extend to higher dimensions with two-way interaction of the gas and rigid body. Another class of methods are the cut-cells methods from \cite{dechriste_numerical_2012, dechriste_cartesian_2016}. These are first-order conservative finite volume type methods for Cartesian grids. Along the moving boundary, grid cells are ``cut'' and the scheme is applied to the now polygon-shaped cell. Additionally, we highlight two finite--volume Arbitrary Lagrangian--Eulerian (ALE) kinetic methods. The first is the Unified Gas--Kinetic Scheme (UGKS) \cite{xu_unified_2010, chen_unified_2012}, in which the kinetic distribution function and its macroscopic moments are evolved simultaneously using the integral solution of the BGK model. The second ALE method \cite{titarevArbitraryLagrangianEulerianDiscrete2022, morozov_dynamics_2022} extends a discrete--velocity finite--volume framework and computes numerical fluxes while treating the collision operator implicitly through a Newton--type iterative procedure. Finally, it should be noted that deterministic methods for rarefied gases are mostly limited to mono-atomic single species cases. When gas mixtures and complicated particle interactions are required, DSMC methods remain very relevant \cite{shrestha_numerical_2015}. Finally, there exists a class of meshfree methods. Meshfree methods don't rely on fixed meshes. Instead, they use an arbitrary (potentially moving) set of points without any information on the connectivity, making them ideal for simulations involving deformable domains, rigid body interactions, and free surfaces.

In this work, we improve upon a meshfree ALE method from \cite{tiwari_meshfree_2022} in several aspects. This method moves spatial points with the mean macroscopic velocity of the gas. As a consequence, the gas-rigid body boundary is tracked exactly, however, the grid becomes irregular. The transport step is then discretised using a meshless Weighted Essentially Non Oscillating (WENO) method from \cite{avesani_new_2014}. The meshless WENO method computes a non-linear weighing of the divergence on several stencils. The BGK-collision operator is treated implicitly like in the semi-Lagrangian methods. The result is a higher-order method that works on arbitrary domains. Recent work on meshfree methods \cite{willems2025b} has shown that the meshless WENO method suffers from stability issues. We replace the meshless WENO method with a meshless MUSCL method that uses only one central stencil. We combine this with the Multidimensional Optimal Order Detection (MOOD) \cite{clain_high-order_2011} method to avoid spurious oscillations at discontinuities. Finally, we use a new boundary routine that does not rely on an iterative procedure or extrapolation. We apply the new method to several test cases in one and two spatial dimensions.

This text is organised as follows. In section \ref{section:model}, we discuss the BGK model for the rarefied gas in Lagrangian formulation and the Newton-Euler equations for rigid body motion. In section \ref{section:schemes}, we discuss the velocity, space and time discretisation. Finally, in section \ref{section:numericalExperiments}, the scheme is applied to several well-known test cases.
 
\section{The model}
\label{section:model}
In section \ref{section:BGKmodel}, we give a summary for the BGK model for rarefied gas. In section \ref{section:ALEFormulation}, we rewrite the BGK equation in the so-called Arbitrary-Lagrangian-Eulerian form. This formulation allows for easy treatment of moving boundaries later on. In this manuscript, we limit ourselves to 1D and 2D simulations. In section \ref{section:ChuReduction}, we use the so-called Chu-reduction to reduce the dimension of the velocity variable to the same dimension of the space variable. This dramatically decreases the computational time. In section \ref{section:rigidBodyModel}, we describe the model for the solid objects and discuss the coupling with the rarefied gas. 

\subsection{The BGK model for a rarefied gas}
\label{section:BGKmodel}

We consider the BGK-Boltzmann equation for rarefied gas dynamics \cite{bhatnagar_model_1954} \begin{align}
	\pdv{f(x, v, t)}{t} + v \cdot \nabla_x f(x, v, t) = \frac{1}{\tau} \left( \maxw - f(x, v, t) \right), \label{eq:BGKBoltzmann}
\end{align} with \(f(x, v, t)\) the population density of the gas as a function of time \( t \in \mathbb{R}^+ \), space \(x \in D\subseteq \mathbb{R}^3 \) and kinetic velocity \(v \in \mathbb{R}^3\). The quantity \( \tau \) is the so-called relaxation time. Compared to the Boltzmann equation, in the BGK equation, the Boltzmann collision operator is replaced by a relaxation operator in which the population density relaxes towards the local Maxwellian distribution \begin{align}
	\maxw = \frac{\rho(x, t)}{ \left( 2\pi R_s T(x, t) \right)^{\frac{3}{2}} } \exp \left(- \frac{ \left| v - U(x, t) \right|^2 }{ 2 R_s T(x, t) } \right). \label{eq:localM}
\end{align} The macroscopic quantities \( \rho(x, t) \in \mathbb{R}, U(x, t) \in \mathbb{R}^3 \) and \(T(x, t) \in \mathbb{R} \) are the density, mean velocity and temperature of the gas. These quantities relate to the first three moments of the population density \(f(x, v, t)\) with respect to \(v\) \begin{align}
	\left( \rho(x, t), \rho(x, t) U(x, t), E(x, t) \right) = \int_{\mathbb{R}^3} \phi(v) f(x, v, t) dv \label{eq:macroscopicQ}
\end{align} 
with \( \phi(v)  = \left( 1, v, {|v|^2}/{2} \right) \). The total energy \(E(x, t)\) is the sum of the internal energy and the kinetic energy \begin{align}
	E(x, t) = \frac{3}{2} \rho(x, t) R_s T(x, t) + \frac{1}{2} \rho(x, t) |U(x, t)|^2.
\end{align} The quantity \(R_s\) is the specific gas constant, which is determined by the universal gas constant and the molar mass of the gas of interest as $R_s = R/m_A$, where $R$ is the universal gas constant, and $m_A$ is the atomic weights, expressed in Kg/mole. The elements in \( \phi(v)\) are referred to as the collision invariants because due to \eqref{eq:localM} and \eqref{eq:macroscopicQ} \begin{align}
	\int_{\mathbb{R}^3} \left( \maxw - f(x, v, t) \right) \phi(v) dv = 0.
\end{align} Consequently, by integration of equation \eqref{eq:BGKBoltzmann} with respect to the collision invariants \( \phi(v) \), a system of conservation laws are obtained \begin{align} \label{eq:BGKConservation}
\begin{split}
	\pdv{\braket{f}}{t} + \nabla_x \cdot \braket{v f} &= 0, \\
	\pdv{\braket{v f}}{t} + \nabla_x \cdot \braket{v \otimes  v f} &= 0, \\
	\pdv{\braket{ \frac{|v|^2}{2} f }}{t} + \nabla_x \cdot \braket{ v\frac{|v|^2}{2} f } &= 0. 
\end{split} \end{align} As the relaxation time \(\tau\) tends to zero, the population density \(f(x, v, t)\) tends to the local Maxwellian \(M\). In this limit, the density \(\rho(x, t)\), mean velocity \(U(x, t)\) and energy \(E(x, t)\) satisfy the compressible Euler equations. The derivation can be formalised using the so-called Chapman-Enskog expansion, in which the population density is written as a Hilbert expansion in the so-called Knudsen number Kn \cite{bardos_fluid_1991}. A zero-th order expansion then yields the compressible Euler equations and a first-order expansion yields the compressible Navier-Stokes equations with the wrong Prandtl number Pr. The ellipsoidal statistical (ES)-BGK model was introduced to solve this problem \cite{holway_kinetic_1965}. Finally, we also note that the BGK and ES-BGK model satisfy the H-theorem, also known as the entropy inequality \cite{andries_gaussian-bgk_2000}. In the remainder of this text, we will focus on the simplest form of the BGK equation \eqref{eq:BGKBoltzmann}. However, we note that the extension to the ES-BGK model and other more complex BGK models like the multi-species BGK model \cite{v_bobylev_general_2018} are straightforward, see for example \cite{boscarino_conservative_2025}. 

\subsection{The ALE formulation}
\label{section:ALEFormulation}
The formulation of the BGK-Boltzmann equation \eqref{eq:BGKBoltzmann} that we will discretise in section \ref{section:schemes} is the so-called arbitrary Lagrangian-Eulerian formulation \cite{tiwari_meshfree_2022} \begin{align}
	\dv{x}{t} &= U_{ALE} \label{eq:langrangianMotion} \\ 
	\dv{f}{t} &= (U_{ALE})-v) \cdot \nabla_x f + \frac{1}{\tau} (\maxw - f), \label{eq:ALEpopulationDynamics}
\end{align} 
where 
\[
\dv{f}{t} := \pdv{f}{t} + U_{ALE} \cdot  \nabla_x f
\]
and 
 \(f = f(x, v, t)\). In equation \eqref{eq:langrangianMotion}, we have introduced motion to the the grid with a velocity \(U_{ALE}\) that is yet to be defined. To compensate for this motion, the term \(U_{ALE} \cdot \nabla_x f\) is added to the population dynamics in equation \eqref{eq:ALEpopulationDynamics}. For all grid points that lie inside the domain, we will set \(U_{ALE}\) to the local mean velocity of the gas \(U(x, t)\). In the case of grid points that lie on the boundary, we distinguish between two cases: points that lie on the outer boundary of $D$, and points that lie the boundary of a rigid object. For the former we set \(U_{ALE} = 0\) assuming the outer boundary is a wall which does not move, and for the latter we set \(U_{ALE} = U_w\), with \(U_w\) the velocity determined by the translation and rotation of the rigid body, see section \ref{section:rigidBodyModel}.

\begin{remark}
	The term \(U_{ALE} \cdot \nabla_x f\) in equation \eqref{eq:ALEpopulationDynamics} is a non-conservative flux. Thus, any numerical scheme for \eqref{eq:ALEpopulationDynamics} cannot conserve \(f(x, v, t)\) and thus also not the density \(\rho(x, t)\), mean velocity \(U(x, t)\) and mean temperature \(T(x, t)\).
\end{remark}

\subsection{Chu-reduced form}
\label{section:ChuReduction}
 We limit ourselves to simulations in one and two spatial dimensions with some symmetry such that we can reduce the dimension of the equation using the so-called Chu-reduction \cite{chu_kinetic-theoretic_1965}. In the one-dimensional case, we consider a slab symmetry. We write the space vector and velocity vector as \(x = (x_1, x_2, x_3)\) and \(v = \left( v_1, v_2, v_3\right)\). We assume that \(v_1\) is aligned with the symmetry axis, and define the velocity \(v_{\perp} = \sqrt{v_2^2 + v_3^2}\). Then, if the initial condition only depends on \(x_1, v_1\) and \(v_{\perp}\), the solution at time \(t_f\) is a function of \(x, v_1\) and \(v_{\perp} \). In this case we can define \begin{align}
	g_1(x, t, v_1) = \int_{\mathbb{R}^2} f(t, x, v)dv_2 dv_3, \; g_2(x, t, v_1) = \int_{\mathbb{R}^2} \left( v_2^2 + v_3^2 \right) f(t, x, v)dv_2 dv_3, \label{eq:chuReduction1D}
\end{align} 
where now by $x$ we mean $x_1$.
Multiplying equation \eqref{eq:ALEpopulationDynamics} by \(1\) and \(v_2^2 + v_3^2\) and integrating with respect to \( \left( v_1, v_2 \right) \in \mathbb{R}^2\), we obtain two equations for \(g_1 = g_1(x, v_1, t)\) and \(g_2 = g_2(x, v_1, t)\) 
\begin{align}
	\dv{g_1}{t} + (v_1 - U(x, t)) \pdv{g_1}{x} = \frac{1}{\tau} \left( G_1 - g_1 \right), \quad \dv{g_2}{t} + (v_1 - U(x, t)) \pdv{g_2}{x} = \frac{1}{\tau} \left( G_2 - g_2\right) \label{eq:ChuALE1D}
\end{align} 
with 
\begin{align}
	G_1 &= \int_{\mathbb{R}^2} \maxw dv_2 dv_3 = \frac{\rho(x, t)}{ \left( 2\pi R_s T(x, t) \right)^{\frac{1}{2}} } \exp \left(- \frac{ \left| v_1 - U(x, t) \right|^2 }{ 2 R_s T(x, t) } \right) \\ 
	G_2 &= \int_{\mathbb{R}^2} \left( v_2^2 + v_3^2 \right) \maxw dv_2 dv_3 = 2 R_s T(x, t) G_1.
\end{align} 
The gain in computational efficiency is clear. If we denote \(N_x\) and \(N_v\) as the amount of degrees of freedom in one space and velocity dimension, then, a 1D simulation using \eqref{eq:BGKBoltzmann} has \(N_v^3 N_x\) points in phase space. When the Chu-reduction is used \eqref{eq:chuReduction1D}, this reduces to \(2 N_x N_v\), thereby significantly decreasing the computational time. The same can be achieved in two space dimensions by assuming 2D geometry, i.e.\ that the system does not depend on one of the coordinates, say $x_3$, and defining \begin{align}
	g_1(x, t, v_1, v_2) = \int_{\mathbb{R}} f(t, x, v) dv_3, \quad g_2(x, t, v_1) = \int_{\mathbb{R}} v_3^2 f(t, x, v_1, v_2) dv_3. \label{eq:chuReduction2D}
\end{align} Then, one obtains 
\begin{align}
	\dv{g_1}{t} + (v - U(x, t)) \nabla_x f = \frac{1}{\tau} \left( G_1 - g_1 \right), \quad \dv{g_2}{t} + (v - U(x, t)) \nabla_x f = \frac{1}{\tau} \left( G_2 - g_2\right) \label{eq:ChuALE2D}
\end{align} with \begin{align}
	G_1 &= \int_{\mathbb{R}^2} \maxw dv_3 = \frac{\rho(x, t)}{2\pi R_s T(x, t)} \exp \left(- \frac{ \left| \begin{bmatrix}v_1 \\ v_2 \end{bmatrix} - U(x, t) \right|^2 }{ 2 R_s T(x, t) } \right) \\ 
	G_2 &= \int_{\mathbb{R}^2} v_3^2 \maxw dv_3 = R_s T(x, t) G_1.
\end{align}

\subsection{Newton-Euler equations for a rigid body}
\label{section:rigidBodyModel}
Consider a rigid body with centre of mass \(X_c \in \mathbb{R}^3\) and mass \(m\) upon which a total force \(\mathcal{F} \in \mathbb{R}^3\) and torque \(\mathcal{T} \in \mathbb{R}^3\) act. Then, the velocity \(V\) of the centre of mass is described by the Newton's second law of motion \begin{align}
	m \dv{V}{t} = \mathcal{F}. 
\end{align} The fundamental law describing the rotation of a rigid body is Euler's second law, or the balance of angular momentum, \begin{align}
	\dv{}{t} L = \dv{}{t}\left(I(t) \omega \right) = \mathcal{T}, \label{eq:momentOfInertia1}
\end{align} where \(L \in \mathbb{R}^3 \) is the angular momentum, \(I(t) \in \mathbb{R}^{3 \cross 3} \) the moment of inertia, and \(\omega \in \mathbb{R}^3 \) the angular velocity, all in a fixed (lab) frame of reference. The moment of inertia can be written as \begin{align}
	I(t) = R(t) I^b R(t)^T,
\end{align} with \(R(t) \in \mathbb{R}^{3 \cross 3} \) the rotation matrix mapping the rigid body coordinate system to lab coordinates. The moment of inertia \(I^b \in \mathbb{R}^3 \) is the (fixed) inertia tensor expressed the rigid body coordinate system. Using these definitions, equation \eqref{eq:momentOfInertia1} can written in the following form, \begin{align}
	I(t) \dv{\omega}{t} + \omega \cross \left(I(t) \omega \right) = \mathcal{T}.
\end{align} The total velocity of a point on the rigid body due to translation and rotation is given by \(U_w = V + \omega \cross (x - X_c)\). The total force \(\mathcal{F}\) and torque \(\mathcal{T}\) are the result of the gas exerting pressure on the surface of the rigid body. We denote \(\psi \in \mathbb{R}^{3 \cross 3}\) the pressure tensor, which is given by \begin{align}
	\psi = \int_{\mathbb{R}^2} (v - U_w) \otimes (v - U_w) f(t, x, v) dv. 
\end{align} The force \(\mathcal{F}\) and torque \(\mathcal{T}\) can be computed from the pressure tensor \(\psi\) \begin{align}
	\mathcal{F} = \int_{\delta S} \left(- \psi \cdot _s \right) dA, \ \mathcal{T} = \int_{\delta S} (x - X_c) \times (-\psi \cdot _s) dA, 
\end{align} with \(_s\) the normal unit vector pointing outward with respect to the rigid body, and \(\delta S\) the surface of the rigid body.

\section{Numerical scheme}
\label{section:schemes}
In this section, we discuss the numerical discretisation method of the Chu-reduced form of the ALE formulation of the BGK-Boltzmann equation. The velocity discretisation method is discussed in section \ref{section:VelocityDiscretization}. In section \ref{section:SpaceDiscretization}, we discuss a \textit{meshless} spatial discretisation method that can deal with the ``Lagrangian'' nature of the grid. We use a splitting method to time integrate the transport step explicitly and the relaxation step implicitly. To avoid spurious oscillations at discontinuities due to the transport step, we use an {\em a posteriori\/} limiting method called Multi dimensional Optimal Order Detection (MOOD). The splitting method and MOOD are discussed in section \ref{section:TimeDiscretization}. In section \ref{section:BoundaryConditions} we discuss the implementation of the diffuse reflective boundary conditions. Finally, in section \ref{section:ImplementationAspects}, we discuss some implementation aspects related to the irregular nature of the spatial grid. 

We use \(i \in\{ 1, \dots, N_x\} \), to denote the index of the grid point at position \(x_i\). The index \(k \in \{ 0, \dots, N_v-1\} \) is used to identify the velocity variable \(v_k\). Finally, the index \(n\) is used for the time variable. We will avoid writing the subscripts ``1'' and ``2'' in \(g_1\) and \(g_2\) since the algorithm for \(g_1\) and \(g_2\) is exactly the same. Using this notation, we write \(g_{i, k}^n \approx g(x_i, v_k, t^n)\).

\subsection{Velocity discretisation: discrete velocity method}
\label{section:VelocityDiscretization}

The velocity variable is discretised according to a simple discrete velocity model (DVM). Consider the one-dimensional Chu-reduced form of the Lagrangian BGK equation \eqref{eq:ChuALE1D}. We truncate the only remaining velocity variable \(v_1\) to the domain \( [-v_{\rm max}, v_{\rm max}] \), which we discretise using \(N_v\) equispaced points. The two equations for \(g_1\) and \(g_2\) become a set of \(2N_v\) equations for \begin{align}
	\dv{g_{i, k}}{t} + (v_k - U_i(t)) \pdv{g_{i, k}}{x} = \frac{1}{\tau} \left( G - g_{i, k} \right)
\end{align} 
where \(v_k = -v_{\rm max} + 2\, k\, v_{\rm max}/(N_v - 1), \; k = 0, \dots, N_v - 1\) and with \(U_i(t)\) an approximation for \(U(x_i, t)\). All integrals of \(v\) can then be approximated using the midpoint rule. In the 2D case \eqref{eq:ChuALE2D}, we use a square velocity grid \( [-v_{\rm max}, v_{\rm max}]^2 \), and discretise velocities \(v_1\) and \(v_2\) in the same way as in the one dimensional case. 

\begin{remark}
	We note that the conservation property \eqref{eq:BGKConservation} and the entropy inequality do not hold for the simple DVM used here \cite{mieussens_discrete_2000}. 
\end{remark}

\begin{remark}
	If \(v_{\rm max}\) is taken large enough such that the distribution functions \(g_{1, k}\) and \(g_{2, k}\) are flat at \(k = 0\) and \(k = N_v-1\), then the midpoint rule is spectrally accurate \cite{trefethen_exponentially_2014}.
\end{remark}

\subsection{Space discretisation: meshless MUSCL}
\label{section:SpaceDiscretization}
Due to the moving nature of the spatial grid, we avoid the use of grid based methods that would require remeshing at every time step. Instead, we discretise the spatial variable using a meshless MUSCL scheme \cite{willems2025b}. This method relies on a Moving Least Squares (MLS) approximation of the spatial derivatives \cite{shepard_two-dimensional_1968, lancaster_surfaces_1981, levin_approximation_1998}. These discrete approximations of the derivative are then substituted into equations \eqref{eq:chuReduction1D} or \eqref{eq:chuReduction2D} to obtain a numerical scheme. We will give the method in one spatial dimension and then briefly discuss the extension to two spatial dimensions.

\begin{remark}
We note that these methods do not require a grid with cells and in theory should work with any set of points covered in the domain. As a result, they are extremely flexible when it comes to complex simulation domains, moving boundaries and free surface flows. However, we emphasize that these spatial discretisations are not conservative. 
\end{remark}

\subsubsection{Moving least squares approximation}
\label{section:MLS}
Consider an irregular set of points \(x_i, i = 1 \dots N \) inside the simulation domain. Each point is assigned a set of neighbours which contains the indices of all the points that lie within a distance \(h_{\rm max}\) from \(x_i\) \begin{align}
	\mathcal{C}_i = \{ j | \Delta x_{ij}| < b_{r} \Delta x = h_{\rm max} \text{ and } i \neq j\},
\end{align} 
with \(\Delta x\) the initial spacing of points and 
\(\Delta x_{ij} = x_j - x_i\). We refer to \(b_{r}\) 
as the ``radius factor''. 
Using the stencil \(\mathcal{C}_i\), 
we compute a local least squares polynomial fit of a test function 
\(u_i \approx u(x_i)\). 
This routine will then be used in section \ref{section:Space1D} to approximate spatial derivatives of \(f(x_i, v, t)\). 
We will give a third-order MLS method in 1D and discuss the generalisation to 2D below. We write down a Taylor polynomial approximation at point \(x_i\) to \(x_j\) 
\begin{align}
	u_j = u_i  + \Delta x_{ij} \dv{u}{x} + \frac{\Delta x_{ij}^2}{2} \dv[2]{u}{x} + \mathcal{O} \left( \Delta x_{ij}^3 \right), \quad \forall j \in \mathcal{C}_i. \label{eq:TaylorExpansion}
\end{align}
Then, by minimizing the \(L^2\) error of the fit, an approximation for the spatial derivatives can be obtained 
\begin{align}
	\left(\dv{\Tilde{u}}{x}, \dv[2]{\Tilde{u}}{x}\right)_{x = x_i} = \argmin_{\alpha,\beta} \sum_j
    w_{ij} \left( u_j - u_i - \Delta x_{ij} \alpha - \frac{\Delta x_{ij}^2}{2} \beta
    \right)^2. \label{eq:MLSMinimisation}
\end{align} 
For stability reasons, we introduce a weight function 
\begin{align}
	w_{ij} = w(x_i, x_j) = \exp \left( - \gamma \Delta x_{ij}^2 \right).
\end{align} We set \(\gamma\) to \({6}/{h_{\rm max}^2}\) so that all points are within the support of the Gaussian. Note that other choices of basis functions and weight functions exist. For an overview of different MLS methods, see \cite{tey_moving_2021, suchde_conservation_2018}. The solution of the minimisation problem \eqref{eq:MLSMinimisation} can be written in the so-called normal equations \begin{align}
	\begin{bmatrix}
		\sum_{j \in \mathcal{C}_i} w_{ij} \Delta x_{ij}^2 & \frac{1}{2} \sum_{j \in \mathcal{C}_i} w_{ij} \Delta x_{ij}^3 \\[2mm]
		\frac{1}{2} \sum_{j \in \mathcal{C}_i} w_{ij} \Delta x_{ij}^3 & \frac{1}{4} \sum_{j \in \mathcal{C}_i} w_{ij} \Delta x_{ij}^4
	\end{bmatrix} \begin{bmatrix}
		\dv{\Tilde{u}}{x} \\[2mm] \dv[2]{\Tilde{u}}{x}
	\end{bmatrix} =  \begin{bmatrix}
		\sum_{j \in \mathcal{C}_i} w_{ij} \Delta x_{ij} \Delta u_{ij} \\[2mm] \frac{1}{2} \sum_{j \in \mathcal{C}_i} w_{ij} \Delta x_{ij}^2 \Delta u_{ij}
	\end{bmatrix}, 
\end{align} with \( \Delta u_{ij} = u_j - u_i \). The solution of the normal equations can then be written explicitly as \begin{align}
	\dv{u}{x} \approx \dv{\Tilde{u}}{x} = \sum_{j \in \mathcal{C}_i} \alpha_{ij} (u_j - u_i), \quad \dv[2]{u}{x} \approx \dv[2]{\Tilde{u}}{x} = \sum_{j \in \mathcal{C}_i} \beta_{ij} (u_j - u_i), \label{eq:MLSCoeff}
\end{align} 
where the coefficients  \(\alpha_{ij}\) and \(\beta_{ij}\) are rational functions of \( \Delta x_{ij}\). We will use the coefficients \(\alpha_{ij}\) and \(\beta_{ij}\) in the meshless MUSCL method to obtain a stable discretisation for the linear transport term in the BGK-Boltzmann equation. In 1D, this MLS method based on a \(k\)-th order Taylor expansion, yields a \(k-\ell\)-th order approximation for the \(\ell\)-th order derivative in space (\(k > \ell\)). The method can be easily extended to higher dimensions in space by using multivariate Taylor expansions.

It should be noted that due to the use of a Taylor expansion in \eqref{eq:TaylorExpansion}, the least squares problem can become ill-conditioned, especially for higher-order approximations. Moreover, by solving the normal equations directly, the condition number of the original least squares problem squares. In practice, the least squares problem should be solved using a QR or an SVD decomposition \cite{trefethen_numerical_1997}. In addition, the distances \(\Delta x_{ij}\) can be scaled with \( {1}/{h_{\rm max}} \), to improve the conditioning further. Using these two improvements, we have achieved fifth-order approximations in 1D and third-order approximations in 2D.

\subsubsection{Meshless MUSCL scheme in one space dimension}
\label{section:Space1D}
In this section, we will use the MLS coefficients to obtain a stable meshless discretisation of the transport term in the Lagrangian form of the BGK-Boltzmann equation. We distinguish between a first-order scheme, that is specifically designed to be positive, and higher-order schemes. We will use the higher-order scheme as a default, and switch to the first-order scheme when spurious oscillations are detected. This is explained in detail in section \ref{section:TimeDiscretization}.

\subsubsubsection{First-order scheme}
\label{section:Space1DFirstOrder}
For the first-order scheme in one space dimension we can make use of a simple generalisation of the classical upwind scheme to the irregular grid case. We define \(\mathcal{U}_{i, k}\) the set of neighbours of point \(i\) that lie in the upwind direction, i.e., \(\mathcal{U}_{i, k} = \{ j \in \mathcal{C}_i | \left(v_k - U(x_i, t)\right) \Delta x_{ij}  < 0 \} \). Using the set \(\mathcal{U}_{i, k}\) and a first-order MLS approximation of the spatial derivative, a simple expression for the transport term in \eqref{eq:ChuALE1D} can be found \begin{align}
	\left(v_k - U(x_i, t)\right) \pdv{g(x_i, v_k, t)}{x} \approx \left(v_k - U_i(t)\right) \frac{\sum_{j \in \mathcal{U}_{i, k}} w_{ij} \Delta x_{ij} \left( g_{j, k}(t) - g_{i, k}(t) \right)}{\sum_{j \in \mathcal{U}_{i, k}} w_{ij} \Delta x_{ij}^2}. \label{eq:Space1DFirstOrder}
\end{align} If we only consider transport and thus ignore the BGK-collision operator, we can obtain a first-order scheme using a forward Euler time integrator \begin{align}
	g_{i, k}^{n+1} = g_{i, k}^{n} \left( 1 + \Delta t (v_k - U_i^n)\frac{\sum_{j \in \mathcal{U}_{i, k}} w_{ij} \Delta x_{ij}}{\sum_{j \in \mathcal{U}_{i, k}} w_{ij} \Delta x_{ij}^2} \right) - \Delta t\left(v_k - U_i^n\right) \frac{\sum_{j \in \mathcal{U}_{i, k}} w_{ij} \Delta x_{ij} g_{j, k}^n}{\sum_{j \in \mathcal{U}_{i, k}} w_{ij} \Delta x_{ij}^2}.
\end{align} Due to the choice of upwind stencil, this scheme is \(L^{\infty}\)-stable under the generalized CFL condition \cite{willems2025b} \begin{align}
	\Delta t \leq \min_{k, i} \frac{\sum_{j \in \mathcal{U}_{i, k}} w_{ij} \Delta x_{ij}^2}{ \left| \left(v_k - U_i^n\right) \sum_{j \in \mathcal{U}_{i, k}} w_{ij} \Delta x_{ij} \right|}. \label{eq:1DStabilityCondition}
\end{align} We note that every stencil has a different configuration of points and therefore every stencil has its own stability condition. The time step is thus determined by the strictest condition. Moreover, because the grid moves, the stability condition changes at every time step. However, for simplicity, we fix the time step using a CFL number for the initial grid. The extension of the proposed scheme to higher order is straightforward and has been adopted in meshless particle methods; see, for example, \cite{wangCompactMovingParticle2023, wangFreeSurfaceBoundary2022}. For the kinetic equations in this paper, however, a higher-order upwind discretisation would require solving a separate least-squares problem for each transport velocity \(v_k\), which may become prohibitively expensive when the number of velocities is large. As an alternative, we therefore adopt a meshfree MUSCL approach based on a central stencil, in which upwinding is introduced via reconstruction. This method is described in detail in the next section.

\begin{remark}
	For a regular grid \(\Delta x_{ij} = \Delta x\) and the weight function set to unity, the CFL condition \eqref{eq:1DStabilityCondition} simplifies to the classical CFL condition.
\end{remark}

\subsubsubsection{Higher-order meshless MUSCL scheme}
\label{section:Space1DHigherOrder}
For the higher-order methods we use the meshless MUSCL scheme that was developed in \cite{willems2025b}. The scheme employs a single central MLS stencil, thereby avoiding the need for multiple separate upwind stencils, such as those used in WENO schemes \cite{avesani_new_2014, avesani_alternative_2021, tiwari_meshfree_2022}. The first-order schemes will act as a fallback method in the case that oscillations are detected. We will illustrate the idea of the meshless MUSCL scheme using as an example a second-order variant. The extension to higher orders and two space dimensions are discussed below. Given the coefficients \(\alpha_{ij}\) and \(\beta_{ij}\) associated to the third-order MLS from section \ref{section:SpaceDiscretization}, we can write the transport term in \eqref{eq:ChuALE1D} as \begin{align}
	\left(v_k - U(x_i, t)\right) \pdv{g(x_i, v_k, t)}{x} \approx \left(v_k - U_i(t)\right) \sum_{j \in \mathcal{C}_i} \alpha_{ij} \left( g_{j, k}(t) - g_{i, k}(t) \right). \label{eq:UnstableCentral}
\end{align} 
The discretisation \eqref{eq:UnstableCentral} corresponds to a central discretisation, which produces an ordinary differential equation with unstable modes \cite{willems2025b}. Instead, we introduce a numerical approximation for the midpoint \(g\left(\frac{x_i + x_j}{2}, j, t\right) \approx g_{i,j, k}(t) \), which we define later, and write the transport term as \begin{align}
	\left(v_k - U(x_i, t)\right) \pdv{g(x_i, v_k, t)}{x} \approx 2\left(v_k - U_i(t)\right)\sum_{j \in \mathcal{C}_i} \alpha_{ij} \left( g_{i,j, k}(t) - g_{i, k}(t) \right). \label{eq:willems2025higherordermeshlessschemeshyperbolicEQ1}
\end{align} Because the stencil ``size'' reduces by a factor two due to the use of the midpoints, a factor two must be introduced in the spatial derivative \eqref{eq:willems2025higherordermeshlessschemeshyperbolicEQ1}. The midpoint \(g_{i,j, k}(t)\) is reconstructed either from point \(i\) or from point \(j\) in an upwind manner \begin{align}
	g\left(\frac{x_i + x_j}{2}, v_k, t\right) \approx g_{ij, k}(t) &= \begin{cases}
		\Bar{g}_{i,j, k}(t), & \text{if } \left(v_k - U(x_i, t)\right) \Delta x_{ij} > 0 \\
		\Bar{g}_{j,i, k}(t), & \text{else }
	\end{cases},
\end{align} with \(\Bar{g}_{ij, k}\) and \(\Bar{g}_{ji, k}\) the reconstructed values (ignoring the time dependence for brevity)\begin{align}
	g\left(\frac{x_i + x_j}{2}, v_k\right) &= g_{i, k} + \frac{\Delta x_{ij}}{2}\pdv{g_{i, k}}{x} + \frac{\Delta x_{ij}^2}{8}\pdv[2]{g_{i, k}}{x} + \mathcal{O}(\Delta x_{ij}^3) \nonumber \\
	\approx \Bar{g}_{ij, k} &= g_{i, k} + \frac{\Delta x_{ij}}{2} \sum_{l \in \mathcal{C}_i} \alpha_{il}(g_{l, k} - g_{i, k}) + \frac{\Delta x_{ij}^2}{8} \sum_{k \in \mathcal{C}_i} \beta_{il}(g_{l, k} - g_{i, k}), \label{eq:1DMUSCLmidpointRec1} \\
	g\left(\frac{x_i + x_j}{2}, v_k\right) &= g_{j, k} - \frac{\Delta x_{ij}}{2}\pdv{g_{j, k}}{x} + \frac{\Delta x_{ij}^2}{8}\pdv[2]{g_{j, k}}{x} + \mathcal{O}(\Delta x_{ij}^3)  \nonumber \\
	\approx \Bar{g}_{ji, k} &= g_{j, k} - \frac{\Delta x_{ij}}{2} \sum_{l \in C_j} \alpha_{jl}(g_{l, k} - g_{j, k}) + \frac{\Delta x_{ij}^2}{8} \sum_{l \in C_j} \beta_{jl}(g_{l, k} - g_{j, k}). \label{eq:1DMUSCLmidpointRec2}
\end{align} Note that the coefficients \(\alpha_{ij}\) and \(\beta_{ij}\) must only be computed once per time step. They can then be reused for every discrete velocity in the computation of the gradient \eqref{eq:willems2025higherordermeshlessschemeshyperbolicEQ1} and the reconstruction of the midpoint \eqref{eq:1DMUSCLmidpointRec1} \eqref{eq:1DMUSCLmidpointRec2}. The extension of this scheme to higher-orders is straightforward. Using a higher-order MLS procedure, a higher-order approximation for the first-order derivatives can be obtained. Then, using all available derivatives from the MLS approximation, the reconstruction can be performed to higher-order. In \cite{willems2025b}, a fourth order scheme was obtained and the stability properties of the meshless MUSCL methods were investigated. Finally, we note that due to the higher-order nature of the scheme, it cannot be TVD stable \cite{leveque_numerical_1992}. In section \ref{section:TimeDiscretization}, we discuss a procedure to deal with this issue.

\subsubsection{Meshless MUSCL scheme in two space dimensions}
\label{section:Space2D}

\subsubsubsection{First-order scheme}
\label{section:Space2DFirstOrder}
\begin{figure}
	\centering

	\begin{tikzpicture}
		% Draw x and y axes centred at the origin
		\draw[->] (-2.5, 0) -- (2.5, 0) node[right] {$x$};
		\draw[->] (0, -2.5) -- (0, 2.5) node[above] {$y$};
		
		% Random points in all quadrants
		\fill (0, 0) circle (2pt) node[below left] {\(i\)};
		\fill (1, 2.3) circle (2pt) node[above right] {\(j\)};
		\fill[blue] (0.5, 1.15) circle (2pt) node[left] {};
		\fill (-2.5, -1.7) circle (2pt);
		\fill (2.2, -1.2) circle (2pt);
		\fill (0.9, -1.5) circle (2pt);
		\fill (2, 1.5) circle (2pt);
		\fill (-2.1, 1) circle (2pt);

		% Draw imaginary cell interface
		\draw[-] (-0.4171, 1.5487) -- (1.4171, 0.7513) node[above] {};
		
		% Draw a vector pointing up and to the right
		\draw[->, thick, red] (0,0) -- (0.3987, 0.9171) node[anchor=east] {$\hat{\eta}_{ij}$};
		\draw[->, thick, red] (0,0) -- (-0.9171, 0.3987) node[anchor=south west] {$\hat{s}_{ij}$};
		% \draw[thick] (-5, 3.75) -- (5, -3.75) node[anchor=north west] {};
		
		% Draw line to node
		% \draw[dashed] (0, 0) -- (2, 1.5);
		
		% Draw the angle 
		\draw [<->] (0:0.7)  arc (0:66.5014:0.7) node [right,pos=0.5] {\(\theta_{ij}\)};
		% \draw[thick] (0,0) arc[start angle=-90, end angle=71.57, radius=1cm];
		% \node at (1.3,0.2) {$\theta_{ij}$};
		
		% Fill region
		% \fill[green!20, opacity=0.5] (-5, 5) -- (-5,-5) -- (0,-5) -- (0, 5) -- cycle;
		% \fill[orange!20, opacity=0.5] (-5, 0) -- (-5,-5) -- (5,-5) -- (5, 0) -- cycle;
	\end{tikzpicture}

	\caption{Illustration of the variables used in the first-order spatial discretisation method.}
	\label{fig:UpwindPraveen}
\end{figure}
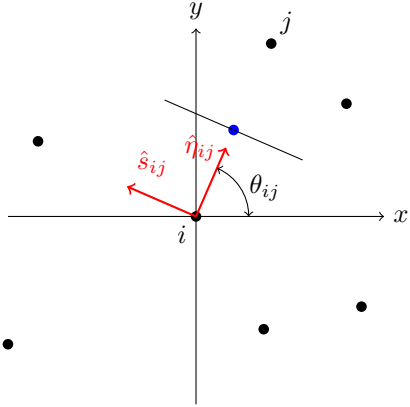

We require a first-order positive scheme to be used as a fallback in case of a MOOD event. Due to potentially negative coefficients resulting from the least squares problem, the first-order scheme from section \ref{section:Space1DFirstOrder} cannot be extended to two spatial dimensions without losing monotonicity. Instead, we use a scheme that was developed in \cite{chandrashekar_positive_2004} and which we summarize here. The scheme is based on a second-order MLS method using a central stencil in two space dimensions. We denote the \(\kappa_{ij}\) and \(\lambda_{ij}\) the MLS coefficients associated to the first order derivatives in the $x$ and $y$ directions, respectively, similarly to what is reported in equation \eqref{eq:MLSCoeff} for the one dimensional case. Then, we define the transport vector \begin{align}
	\vec{a}_i = \begin{bmatrix}
		a_{1, k} \\  a_{2, l}
	\end{bmatrix} = \begin{bmatrix}
		v_{1, k} \\  v_{2, l}
	\end{bmatrix} - U_i(t) = \begin{bmatrix}
	v_{1, k} \\  v_{2, l}
	\end{bmatrix} - \begin{bmatrix}
		U_{1, i}(t) \\ U_{2, i}(t)
	\end{bmatrix}, 
\end{align} discretised according to section \ref{section:VelocityDiscretization}. Next, we define the angle between the line connecting the centre point \(i\) and point \(j\) to be \(\theta_{ij}\), and the unit normal vectors \(\hat{\eta}_{ij} = \left( \cos \theta_{ij}, \sin \theta_{ij} \right)^T \) and \(\hat{s}_{ij} = \left( - \sin \theta_{ij}, \cos \theta_{ij} \right)^T\). The direction \(\hat{\eta}_{ij}\) can be interpreted as the direction orthogonal to an imaginary cell boundary separating point \(i\) and point \(j\), see figure \ref{fig:UpwindPraveen}. The scheme then approximates the flux in the direction \(\hat{\eta}_{ij}\) using an upwind flux. The flux along the imaginary cell boundary orthogonal to \(\hat{s}_{ij}\) is the average flux with some additional diffusion such that the scheme is \(L^{\infty}\)-stable under a generalised time step restriction. Ultimately, the scheme can be written as \begin{align}
	 \vec{a}_i \cdot \nabla_x g(x_i, v, t)  &\approx \sum_{j \in \mathcal{C}_i} \left[ \Bar{\kappa}_{ij} \left( \vec{a}_i \cdot \hat{\eta}_{ij} - |\vec{a}_i \cdot \hat {\eta}_{ij}| \right) + \left( \Bar{\lambda}_{ij} (\vec{a}_i \cdot \hat{s}_{ij}) - \left| \Bar{\lambda}_{ij} (\vec{a}_i \cdot \hat{s}_{ij}) \right| \right)  \right]\left(g_{j, kl} - g_{i, kl}\right) \\ &=: -\sum_{j \in \mathcal{C}_i} c_{ij} \left(g_{j, kl} - g_{i, kl}\right),  \label{eq:2DFirstOrderTransport}
\end{align}
with the first set of indices of \(g_{j, kl}\) referring to the position and the second set of indices referring to the velocity and 
\begin{align}
	\begin{bmatrix}
		\Bar{\kappa}_{ij} \\ \Bar{\lambda}_{ij}
	\end{bmatrix} = \begin{bmatrix}
		\hat{\eta}_{ij} & \hat{s}_{ij}
	\end{bmatrix}^T \begin{bmatrix}
		\kappa_{ij} \\ \lambda_{ij}
	\end{bmatrix}.
\end{align} 
The coefficients \(\Bar{\kappa}_{ij}\) and \(\hat{\eta}_{ij}\) provide an approximation of the derivatives along the direction \(\hat{\eta}_{ij}\) and \(\hat{s}_{ij}\) respectively. A first-order in time scheme for the transport in \eqref{eq:2DFirstOrderTransport} is \(L^{\infty}\)-stable provided that \begin{align}
	\Delta t \leq \min_i \frac{1}{\sum_{j \in \mathcal{C}_i} c_{ij}}.
\end{align}
The extension of the scheme to three space dimensions is discussed in \cite{mariappan_high-performance_2024}. 

\subsubsubsection{Second-order scheme}
\label{section:Space2DSecondOrder}
For the higher-order scheme in 2D we use a straightforward extension of the meshless MUSCL scheme from section \ref{section:Space1DHigherOrder}. As an example, we will summarize a second-order method here. Consider a 2D third-order MLS approximation analogous to section \ref{section:MLS}, that yields the following local approximations of the spatial derivatives \begin{align}
	\pdv{g_{i, kl}}{x} &\approx \pdv{\Tilde{g}_{i, kl}}{x} = \sum_{j \in \mathcal{C}_i} \xi_{ij} (g_{j, kl} - g_{i, kl}), \; \pdv{u_i}{y} \approx \pdv{\Tilde{g}_{i, kl}}{y} = \sum_{j \in \mathcal{C}_i} \zeta_{ij} (g_{j, kl} - g_{i, kl}), \label{eq:approx1} \\ 
	\pdv{g_{i, kl}}{x}{y} &\approx \pdv{\Tilde{g}_{i, kl}}{x}{y} = \sum_{j \in \mathcal{C}_i} \gamma_{ij} (g_{j, kl}  - g_{i, kl}), \; \pdv[2]{u_i}{x} \approx \pdv[2]{\Tilde{g}_{i, kl}}{x} = \sum_{j \in \mathcal{C}_i} \eta_{ij} (g_{j, kl} - g_{i, kl}), \\
	\pdv[2]{g_{i, kl}}{y} &\approx \pdv[2]{\Tilde{g}_{i, kl}}{y} = \sum_{j \in \mathcal{C}_i} \nu_{ij} (g_{j, kl}  - g_{i, kl}).
\end{align} We can then write the divergence in \eqref{eq:ChuALE2D} as \begin{align}
	a_{1, k} \pdv{g(x_i, v_k, t)}{x} + a_{2, l} \pdv{g(x_i, v_k, t)}{y} \approx 2 a_{1, k} \sum_{j \in \mathcal{C}_i} \xi_{ij} \left(g_{ij, kl} - g_{i, kl}\right) + 2 a_{2, l} \sum_{j \in \mathcal{C}_i} \zeta_{ij} \left(g_{ij, kl} - g_{i, kl}\right), 
\end{align} with the midpoint \(g_{ij, kl}\) is defined as \begin{align}
	g_{ij, kl} &= \begin{cases}
		\Bar{g}_{ij, kl}, & \text{if } a_x \Delta x_{ij} + a_y \Delta y_{ij} > 0 \\
		\Bar{g}_{ji, kl}, & \text{otherwise. } 
	\end{cases}
\end{align} The reconstructions from the centre point \(\Bar{g}_{ij, kl}\) and the neighbouring point \(\Bar{g}_{ji, kl}\) are computed using a Taylor polynomial in which we have replaced the exact derivatives with their MLS approximations \begin{align}
	\Bar{g}_{ij, kl} &= g_{i, kl} + \frac{\Delta x_{ij}}{2}\pdv{\Tilde{g}_{i, kl}}{x} + \frac{\Delta y_{ij}}{2}\pdv{\Tilde{g}_{i, kl}}{y} + \frac{\Delta x_{ij}^2}{8}\pdv[2]{\Tilde{g}_{i, kl}}{x} + \frac{\Delta y_{ij}^2}{8}\pdv[2]{\Tilde{g}_{i, kl}}{y} + \frac{\Delta x_{ij} \Delta y_{ij}}{4} \pdv{\Tilde{g}_{i, kl}}{x}{y}, \label{eq:2DMUSCLmidpointRec1} \\
	\Bar{g}_{ji, kl} &= g_{j, kl} - \frac{\Delta x_{ij}}{2}\pdv{\Tilde{g}_{j, kl}}{x} - \frac{\Delta y_{ij}}{2}\pdv{\Tilde{g}_{j, kl}}{y} + \frac{\Delta x_{ij}^2}{8}\pdv[2]{\Tilde{g}_{j, kl}}{x} + \frac{\Delta y_{ij}^2}{8}\pdv[2]{\Tilde{g}_{j, kl}}{y} + \frac{\Delta x_{ij} \Delta y_{ij}}{4} \pdv{\Tilde{g}_{j, kl}}{x}{y} . \label{eq:2DMUSCLmidpointRec2}
\end{align} The MLS derivatives are second-order approximations. The midpoints used in the MLS derivatives are third-order approximations. The full scheme is thus of second-order. 

\begin{remark}
	We note that meshfree schemes that rely directly on least-squares polynomial reconstructions to discretise differential operators are inherently non-conservative. Nevertheless, recent work on meshfree methods \cite{willems2025b} has demonstrated that the meshless MUSCL approach adopted here substantially reduces non-conservation errors compared to existing techniques, including the meshless WENO scheme of \cite{tiwari_meshfree_2022}. In Section~\ref{section:shockTube}, we present several numerical experiments that examine and quantify the lack of conservation.
\end{remark}

\subsection{Time discretisation: IMEX and MOOD} % IMEX, MOOD, timestep restriction
\label{section:TimeDiscretization}
We discretise time using 
%Diagonally Implicit 
IMEX Runge-Kutta schemes 
%(DIRK) 
in which the transport term is discretised explicitly and the BGK collision term is discretised implicitly by a Diagonally Implicit (DIRK) scheme. We denote \(\tilde{a}_{sp}\) and \(\tilde{w}_{s}\) the coefficients associated to the explicit Runge-Kutta scheme, and \(a_{sp}\) and \(w_{s}\), the coefficients associated to the implicit Runge-Kutta scheme. The equation that describes the motion of the particles \eqref{eq:langrangianMotion} is discretised using the explicit scheme of the IMEX Runge-Kutta method. Using the notation \(g(t^n, x_i, v_{k}) \approx g^n_{i, kl}\), we obtain the following fully discrete scheme 
\begin{align}
    \Tilde{g}^{(s)}_{i, k} &= g^n_{i, k} - \Delta t \sum_{p=1}^{s-1} \tilde{a}_{sp} \transp (g_{i, k}^{(p)})\label{eq:IMEX4}\\
	\bar{g}^{(s)}_{i, k} &= \Tilde{g}^{(s)}_{i, k} + \frac{\Delta t}{\tau} \sum_{p=1}^{s-1} a_{sp} \left( G^{p}_{i,k} - g^{(p)}_{i,k} \right), 
    \quad x^{(s)} = x^n + \Delta t \sum_{p=1}^{s-1} \tilde{a}_{sp} U_{i}^p,\label{eq:IMEX1}\\
    \quad g^{(s)}_{i, k} & = \bar{g}^{(s)}_{i, k} + \frac{\Delta t}{\tau}  a_{ss} \left( G^{s}_{i, k} - g^{(s)}_{i, k} \right), \label{eq:IMEX2}\\
	g^{n+1}_{i, k} &= g^n_{i, k} - \Delta t \sum_{s=1}^{\nu} \tilde{w}_{s} \transp (g_{i, k}^{(s)}) + \frac{\Delta t}{\tau} \sum_{s=1}^{\nu} w_{s} \left( G^{s}_{i,k} - g^{(s)}_{i,k} \right), 
    \quad x^{n+1} = x^n + \Delta t \sum_{s=1}^{\nu} \tilde{w}_{s} U_{i}^s, \label{eq:IMEX3}
\end{align} with \(\transp\) a meshless MUSCL discretisation of the transport operator. Notice that, because of the DIRK structure of the implicit scheme, all terms in Eq.\eqref{eq:IMEX1} can be computed explicitly. The implicit step \eqref{eq:IMEX2} can be computed in an explicit fashion taking advantage of the special structure of the BGK model, by using a technique that was first published in \cite{pieraccini_implicitexplicit_2007} and is often used for discrete velocity models of the BGK equation \cite{cho_conservative_2022, carrillo_conservative_2024}. The idea behind the technique is very simple. We multiply \eqref{eq:IMEX2} by $(1,v_k,|v_k|^2)$ and sum over $k$ \begin{align}
	\sum_k g^{(s)}_{i, k} (1,v_k,|v_k|^2) = \sum_k \left(\bar{g}^{(s)}_{i, k} + \frac{\Delta t}{\tau}  a_{ss} \left( G^{s}_{i, k} - g^{(s)}_{i, k} \right)\right) (1,v_k,|v_k|^2). \label{eq:BGKIMEXTrick}
\end{align} By definition of local Maxwellian, \begin{align*}
	\sum_k(G^s_{i,k}-g^{(s)}_{i,k})(1,v_k,|v_k|^2) = 0
\end{align*} and therefore the moments of $g^{(s)}_{i,k}$ in \eqref{eq:BGKIMEXTrick} are computed explicitly as moments of $\bar{g}^{(s)}_{i,k}$. Once the moments are computed, we can explicitly compute the Maxwellian $G_{i,k}^s$, and finally 
\[
    g^{(s)}_{i,k} = \frac{\tau \bar{g}^{(s)}_{i,k} + \Delta t a_{ss} G^s_{i,k}}{\tau+\Delta t \, a_{ss}}.
\]

For the numerical test in section \ref{section:numericalExperiments}, we will use the ARS(2, 2, 2) scheme \cite{ascher_implicit-explicit_1997} and the IMEX-SSP2(3, 3, 2) scheme \cite{pareschi_implicitexplicit_2005}.

To avoid spurious oscillations at discontinuities due to the discretisation of the transport term, we employ the Multi-dimensional Optimal Order Detection algorithm, commonly referred to as MOOD \cite{clain_high-order_2011}. MOOD uses an \emph{a posteriori} check to detect oscillations. Upon detection of an oscillation, the order of accuracy of the spatial discretisation method is reduced to one. The criterion we use to detect oscillations is the discrete maximum property (DMP) \cite{bressan_hyperbolic_2005}. This is a natural choice since the transport in the BGK equation is scalar and linear, thus, for the solution to be a weak solution, it must satisfy the DMP.  Numerically, we check if the solution after transport \(\Tilde{g}_{i, k}^{(s)}\) \eqref{eq:IMEX4} does not exceed the minimum and maximum values from the previous stage 
\begin{align}
	\min_{j \in \mathcal{C}_i} (g_{i, k}^{(s-1)}, g_{j, k}^{(s-1)}) \leq \Tilde{g}_{i, k}^{(s)} \leq \max_{j \in \mathcal{C}_i} (g_{i, k}^{(s-1)}, g_{j, k}^{(s-1)}). \label{eq:MOODDMP}
\end{align} If this criterion is not satisfied, the order of the method is reduced, and the solution \(\Tilde{g}_{i, k}^{(s)}\) is recomputed. However, in \cite{diot_improved_2012}, it was shown that the DMP is too strict at smooth extrema, thereby needlessly reducing the order of the method. We solve this problem by supplementing the DMP condition with the so-called u2 detection criterion. This criterion relaxes the DMP at local extrema to prevent losing the order of the method. First, we define the following smoothness indicators \begin{align}
\Tilde{\mathcal{X}}_{i, k}^{\rm min} &= \min_{j \in \mathcal{C}_i} \left( \left| \pdv[2]{\Tilde{g}_{i, k}}{x} \right| , \left| \pdv[2]{\Tilde{g}_{j, k}}{x} \right| \right) \text{ and } \Tilde{\mathcal{X}}_{i, k}^{\rm max} = \max_{j \in \mathcal{C}_i} \left( \left| \pdv[2]{\Tilde{g}_{i, k}}{x} \right|, \left| \pdv[2]{\Tilde{g}_{j, k}}{x} \right| \right), \label{eq:oscillationX1} \\
\mathcal{X}_{i, k}^{\rm min} &= \min_{j \in \mathcal{C}_i} \left( \pdv[2]{\Tilde{g}_{i, k}}{x}, \pdv[2]{\Tilde{g}_{j, k}}{x} \right) \text{ and } \; \mathcal{X}_{i, k}^{\rm max} = \max_{j \in \mathcal{C}_i} \left( \pdv[2]{\Tilde{g}_{i, k}}{x}, \pdv[2]{\Tilde{g}_{j, k}}{x} \right). \label{eq:oscillationX2} \end{align} The smoothness indicators are based on the minimal and maximal values of the curvatures in a stencil. Then, the u2 detection criterion consists of three conditions \begin{align}
	\left| \max_{j \in \mathcal{C}_i} (g_{i, k}^{(p-1)}, g_{j, k}^{(p-1)}) - \min_{j \in \mathcal{C}_i} (g_{i, k}^{(p-1)}, g_{j, k}^{(p-1)}) \right| &< \delta \label{eq:MOODFlat}\\
	\mathcal{X}_{i, k}^{\rm min} \mathcal{X}_{i, k}^{\rm max} &> - \delta \label{eq:MOODSignCheck} \\
	\frac{\tilde{\mathcal{X}}_{i, k}^{\rm min}}{\tilde{\mathcal{X}}_{i, k}^{\rm max}} \geq 1/2 \text{ or } \tilde{\mathcal{X}}_{i, k}^{\rm max} &< \delta \label{eq:MOODOscillationCheck}
\end{align} The first condition \eqref{eq:MOODFlat} prevents MOOD events in regions where the solution is flat. Criterion \eqref{eq:MOODSignCheck} checks if a point is at a local extremum by checking the signs of the minimum and maximum curvatures. Finally, criterion \eqref{eq:MOODOscillationCheck} is used to distinguish between discontinuities and smooth extrema. All criteria contain an additional parameter \(\delta\) to avoid micro-oscillations triggering MOOD events and thus resulting in unnecessary loss of accuracy. In two space dimensions, the criteria \eqref{eq:MOODOscillationCheck} and \eqref{eq:MOODSignCheck} are checked in both space dimensions, see appendix \ref{appendix:MOOD2D}.

The full MOOD procedure can be summarised as follows. First, the transport step \eqref{eq:IMEX4} is computed using a higher-order method. Then, the discrete maximum property \eqref{eq:MOODDMP} is checked. If this criterion is satisfied, the implicit part of the IMEX stage is computed according to  (\ref{eq:IMEX1}-\ref{eq:IMEX3}). If the discrete maximum property is not satisfied, the u2 conditions are checked according to (\ref{eq:MOODFlat}-\ref{eq:MOODOscillationCheck}). If all three conditions are satisfied, the solution is still accepted, and the implicit part of the IMEX scheme is computed. If one of the u2 criteria fails, the solution \(\Tilde{g}^{(s)}_{i, k}\) is recomputed using a forward Euler time integration method and a first-order method in space \begin{align}
    \Tilde{g}^{(s)}_{i, k} &=  g^{s-1}_{i, k} - \Delta t \left(w_s - w_{s-1}\right) \transp(g_{i, k}^{(s-1)}).
\end{align} Now, the remainding computations in the IMEX stage can be computed using again  (\ref{eq:IMEX1}-\ref{eq:IMEX3}). 
\begin{remark}	
	We note that it is possible to use a Strong Stability Preserving (SSP) IMEX Runge-Kutta method such that it is not necessary to switch to explicit Euler time integration when a MOOD event occurs. We did not achieve as good results using an SSP IMEX method versus switching to forward Euler. 
\end{remark}

\subsection{Diffuse-reflective boundary conditions}
\label{section:BoundaryConditions}
In all the numerical examples we will consider diffuse reflective boundary conditions. Let us denote by 
$\partial D$ the boundary of the domain occupied by the gas, and \(\hat n\) the unit normal to the boundary pointing inside the gas domain. At the particle level, diffusive reflective boundary conditions imply that upon collision with a wall, the particle is reflected back into the gas with a velocity sampled from a Maxwellian with mean velocity \(U_w\) and temperature \(T_w\). The boundary condition reads \begin{align}
	f(x, v, t) = \frac{\rho_w}{ \left( 2 \pi R_s T_w \right)^{\frac{3}{2}} } \exp \left( -\frac{|v - U_w|^2}{2 R_s T_w} \right), \quad x \in \partial D, \; (v - U) \cdot \hat n \geq 0. \label{eq:DRBC}
\end{align} In the case of a moving object, the velocity of the wall \(U_w\) is determined based on the translation and rotation of the object, see section \ref{section:rigidBodyModel}. The density of the gas \(\rho_w\) that is reflected back into the domain is determined by mass conservation at the surface of the wall 
\begin{align}
	 \int_{(v - U_w) \cdot \hat n > 0} M_w (v-U_w) \cdot \hat n \; dv + \int_{(v - U_w) \cdot \hat n < 0} f(x, v, t) (v-U_w) \cdot \hat n \; dv = 0. \label{eq:DRBCZeroMassFlux}
\end{align} 
The first term in equation \eqref{eq:DRBCZeroMassFlux} represents the particles that have been reflected off the wall, the second represents the particles that are headed towards the wall. In \cite{groppi_boundary_2016}, an algorithm is discussed for numerically imposing this boundary condition in the context of semi-Lagrangian methods. They enforce condition \eqref{eq:DRBC} after the collision step, which ultimately requires an iterative procedure to be used. Alternatively, they present an algorithm in which the distribution function is extrapolated from points inside the domain to grid points on the boundary, after which the ``wall density'' \(\rho_w\) \eqref{eq:DRBCZeroMassFlux} can be computed. We perceive both options to be suboptimal: the iterative procedure is cumbersome and can be computationally expensive (especially for small Knudsen numbers), and the polynomial extrapolation procedure is potentially unstable. We therefore propose an alternative procedure to implement diffuse reflective boundary conditions that neither requires an extrapolation procedure nor an iterative procedure. We propose to enforce the boundary condition \eqref{eq:DRBC} after a transport step and thus before the collision step. After all, gas particles interact with the wall during the transport step. The algorithm consists of three steps that are applied in each Runge-Kutta stage to points that lie on the boundary of the gas domain. For simplicity, we give the numerical procedure using a first-order time integration algorithm and we leave the space and velocity variable continuous. First, we compute the transport step for all ``outgoing'' characteristics
\begin{align}
	\Bar{f} = f^n - \Delta t(v - U_w) \cdot  \nabla_x f^n  \quad \text{for} \quad  (v - U_w) \cdot \hat n < 0 \label{eq:DRBCStep1},
\end{align} where, as in section \ref{section:TimeDiscretization}, the quantity \( \Bar{f} \) is the distribution function after advection at time \(n+1\), but before collisions have been taken into account. The density $\rho_w$ that, together with the boundary temperature $T_w$ and velocity $U_w$, characterizes the Maxwellian distribution of particles reflected back into the gas domain, can be computed imposing that the net mass flux across the boundary is zero $\forall x\in\partial D$,   using \eqref{eq:DRBCZeroMassFlux}. Hence, numerically this becomes \begin{align}
	\int_{(v - U_w) \cdot \hat n > 0} (v-U_w) \cdot \hat n \; M^{n+1}(\rho_w, U_w, T_w) dv + \int_{(v - U_w) \cdot \hat n < 0} (v-U_w) \cdot \hat n \; \Bar{f} dv = 0
\end{align} 
Finally, we assign the ``incoming'' half of the distribution function as
\begin{align}
	\Bar{f} = M(\rho_w, U_w, T_w) \quad \text{for} \quad  (v - U_w) \cdot \hat n > 0, 
\end{align} 
after which the implicit step of the Runge-Kutta scheme can be applied to both the grid points inside the domain and on the boundary. This boundary procedure is not limited to meshless discretisations and could also be applied to semi-Lagrangian schemes.

\subsection{Implementation aspects}
\label{section:ImplementationAspects}
An important aspect of meshless methods is the so-called grid management\footnote{Although the methods adopted in this paper are denoted as ``meshless" or ``gridless", we still denote by ``grid" the set of points, to emphasize the idea that the smoothness of their distribution is relevant for the spacial accuracy of the method. For this reason we still use the term ``grid". } \cite{suchde_conservation_2018, seibold_m-matrices_2006, tiwari_meshfree_2022}. Grid management involves creating the stencils \(\mathcal{U}_{i, k}\) and \(\mathcal{C}_{i}\). Moreover, the stencils \(\mathcal{U}_{i, k}\) and \(\mathcal{C}_{i}\) must be recomputed regularly since the grid moves. If these stencils are not updated regularly enough, the spatial discretisation errors may become very large. A grid management routine also regularly adds or removes points from the domain to avoid that the grid becomes too irregular. For example, in case there is part of the domain that is void of grid points, a new point is inserted. Similarly, when multiple points are too close together, they are replaced by a point in the middle. As was investigated in \cite{willems2025b}, this improves the stability of the schemes. When a new point is added or removed, the stencils \(\mathcal{U}_{i, k}\) and \(\mathcal{C}_{i}\) of all points should be recomputed. 

\begin{figure}	
	\centering
	\hspace*{\fill}  % Horizontally center
	\begin{tikzpicture}[scale=0.8]
		
		% Set the seed for reproducibility
		\pgfmathsetseed{2023}
		
		% Grid parameters
		\def\gridSize{8}
		\def\spacing{1}
		
		% Highlight one cell (e.g., cell at i=3, j=2)
		\fill[green!30] (2*\spacing, 2*\spacing) rectangle (7*\spacing, 7*\spacing);
		
		% Draw vertical lines
		\foreach \x in {0,...,\gridSize} {
			\draw[gray!70] (\x*\spacing, 0) -- (\x*\spacing, \gridSize*\spacing);
		}
		
		% Draw horizontal lines
		\foreach \y in {0,...,\gridSize} {
			\draw[gray!70] (0, \y*\spacing) -- (\gridSize*\spacing, \y*\spacing);
		}
		
		% Annotate spacing
		\draw[<->, thick] (0, -0.5) -- (\spacing, -0.5);
		\node at (0.5*\spacing, -0.9) {$b_{v} \Delta x$};
		
		% Place one random point inside each cell
		\foreach \i in {0,...,7} {
			\foreach \j in {0,...,7} {
				\pgfmathsetmacro{\dx}{rnd}
				\pgfmathsetmacro{\dy}{rnd}
				\fill[blue] ({(\i+\dx*0.95)*\spacing}, {(\j+\dy*0.95)*\spacing}) circle (2pt);
			}
		}
		
		% Randomly placed 16 points on grid intersections
		\foreach \i in {1,...,10} {
			\pgfmathsetmacro\x{(rand+1)*4)}
			\pgfmathsetmacro\y{(rand+1)*4)}
			\fill[blue] (\x*\spacing, \y*\spacing) circle (2pt);
		}
		
		% Fixed point			
		\pgfmathsetmacro{\x}{4.8}
		\pgfmathsetmacro{\y}{4.2}
		\def\r{2.1}
		\def\angle{-15}
		
		\fill[red] (\x, \y) circle (2pt);
		\draw (\x, \y) circle (\r);
		
		\draw[<->, thick] (\x,\y) -- ++(\angle:\r);
		\node[anchor=east, rotate=\angle] at ($(\x,\y)!0.5!(\x,\y)+(\angle:\r) + (0, 0.3)$) {$b_{r} \Delta x$};

	\end{tikzpicture}
	\hspace*{\fill}  % Horizontally center
	\caption{To find the neighbours of the red grid point, the distance to all points in the green highlighted cells are checked. Only points inside the circle are assigned a neighbour to the red point. }
	\label{fig:particleMangement}
\end{figure}
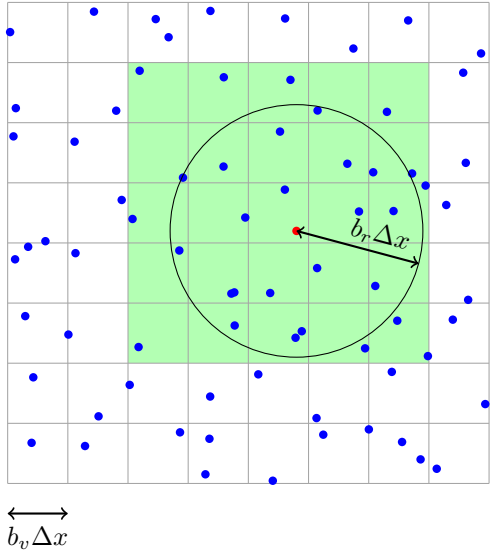

In one spatial dimension, the grid points are stored in an array that is sorted using the position. In this way, it is trivial to efficiently find neighbouring points that are close, and to add and remove points. In two spatial dimensions the grid points are stored in an array in arbitrary order. To find neighbouring points, we create a uniform background square grid with spacing \(b_{v} \Delta x\), \(b_v \geq 1\), with \(\Delta x\) the initial grid spacing. We refer to \(b_{v}\) as the ``voxel factor''. Using the uniform background grid, we generate a table that contains the indices of the grid points per background cell. When the stencil \(\mathcal{C}_{i}\) must be recomputed, we search \(\ceil*{{b_r}/{b_v}}\) cells in all directions for neighbouring points (see figure \ref{fig:particleMangement}), thereby avoiding a costly \(N_x^2\) search. Using the neighbour information of every point, two neighbouring points are removed and replaced with a point in the middle if the distance between the points is less than \(b_{\rm minDist} \Delta x\) with \(b_{\rm minDist} < 1\). Finally, when a background cell is void of grid points, a new grid point is placed at its cell centre, to avoid gaps in the domain.

When a new point is added to the domain, the solution value at that point is interpolated using a second-order MLS interpolation procedure. When MOOD is used, the discrete maximum property \eqref{eq:MOODDMP} is checked for the newly interpolated value. In case the DMP criterion is not satisfied, we recompute the solution value using a first-order interpolation routine.
\begin{remark}
	Due to the adaptive nature of the grid, throughout the simulation, the grid can contain anywhere between 
    \({N_x^0}/{b_v}\) and \({N_x^0}/{b_{\rm minDist}}\) grid points in one spatial dimension with \(N_x^0\) the initial amount of grid points.
\end{remark}

\section{Numerical experiments}
\label{section:numericalExperiments}
In this section, we perform several simulations with the newly developed numerical schemes. To demonstrate the correctness of the algorithms, we do a convergence analysis in one and two spatial dimensions. The behaviour of the scheme around shocks is illustrated using Sod's shock tube example. To show the versatility of the scheme when it comes to moving boundaries, we consider two simulations: the one-dimensional moving plate problem and the driven square cavity with rigid object \cite{tiwari_meshfree_2022}. Due to the moving nature of the grid, simulations with moving boundaries can be tackled directly, and the boundary of the gaseous domain is tracked exactly with no need for interpolation. Finally, we consider the shear layer problem from \cite{banda_lattice-boltzmann_2007, kurganov_third-order_2000, bell_second-order_1989}. The code, along with all examples, is available at \cite{githubRepo1} and \cite{githubRepo2}. 
For all simulations, the grid is uniform at the initial time, and the simulation parameters are given in table \ref{tab:globalParams}. 

We compare the meshless MUSCL methods with a second-order WENO MLS method from \cite{tiwari_meshfree_2022}. In 1D, the WENO method uses a non-linear combination of three second-order polynomials: two one-sided polynomials (left and right) and one central polynomial. The non-linear weights are chosen as in \cite{avesani_new_2014}. MOOD is not activated when WENO is used. Compared to the meshless WENO scheme, the MUSCL method incurs an additional cost due to reconstruction in the midpoints. In one dimension, this cost is non-negligible. However, in two space dimensions, the meshless WENO scheme requires multiple stencils in each dimension, making the meshless MUSCL scheme more efficient. We refer to \cite{willems2025b} for a detailed comparison between the two meshless schemes for the linear advection case. 

\begin{table}[h!]
	\centering
	\begin{tabular}{|c||c|c|c|c|} 
		\hline
		& $b_{\rm minDist}$ & $b_v$ & $b_r$ & \(\delta \) \\ \hline \hline
		1D & 0.1 & 2.0 & 4.0 & $1\times10^{-7}$ \\
		2D & 0.5 & 1.0 & 2.5 & $1\times10^{-8}$ \\
		\hline
	\end{tabular}
	\caption{The parameter \(b_{\rm minDist} \Delta x\) determines the minimum distance that is maintained between two points. The parameter $b_v \Delta x$ in one spatial dimension is the maximum distance between two points. In two spatial dimensions it is the size of the background grid. The parameter \(b_r \Delta x\) determines the size of the stencil \(\mathcal{C}_i\). Finally, \(\delta\) is a fixed parameter used in MOOD to avoid micro oscillations, see section \ref{section:TimeDiscretization}.}
	\label{tab:globalParams}
\end{table}

\subsection{Convergence in one space dimension}
In this section, we perform a convergence study of the meshless schemes described in section \ref{section:schemes}. Consider a fixed one-dimensional domain \([-1, 1]\), inside which there is a gas with gas constant \(R_s\) equal to one and relaxation time set to \(1 \times 10^{-5}\). The gas is initially in local thermodynamic equilibrium with mean velocity 
\begin{align}
	U(x) = \frac{e^{-5(x - 0.1)^2} - 2 e^{-5(x+0.3)^2}}{10},
\end{align} 
and the temperature and density are set equal to one. We apply diffuse-reflective boundary conditions at the fixed boundaries with temperature set to one. The velocity variable is discretised using \(N_v = 20\) equidistant points in the domain \([-10, 10]\). The simulation is performed up to time \(t_f=0.04\) (before the onset of shocks) for a range of spatial discretisation points \(N_x\). We set the CFL number to \(0.02\) such that the time discretisation error is sufficiently small. We compute the relative \(L^1\) error in the moments, one for every combination of algorithm and grid by interpolating the numerical solution from the irregular grid to a uniformly distributed fixed grid with 500 grid points using fifth-order spline interpolation. The solution with which the error is computed is generated using the IMEX-SSP2(3, 3, 2) scheme with the fourth order meshless MUSCL method using \(1500\) grid points. 

\begin{figure}[h!]
	\centering
	\begin{subfigure}[b]{0.49\textwidth}
		\centering
		\includegraphics[width=\linewidth]{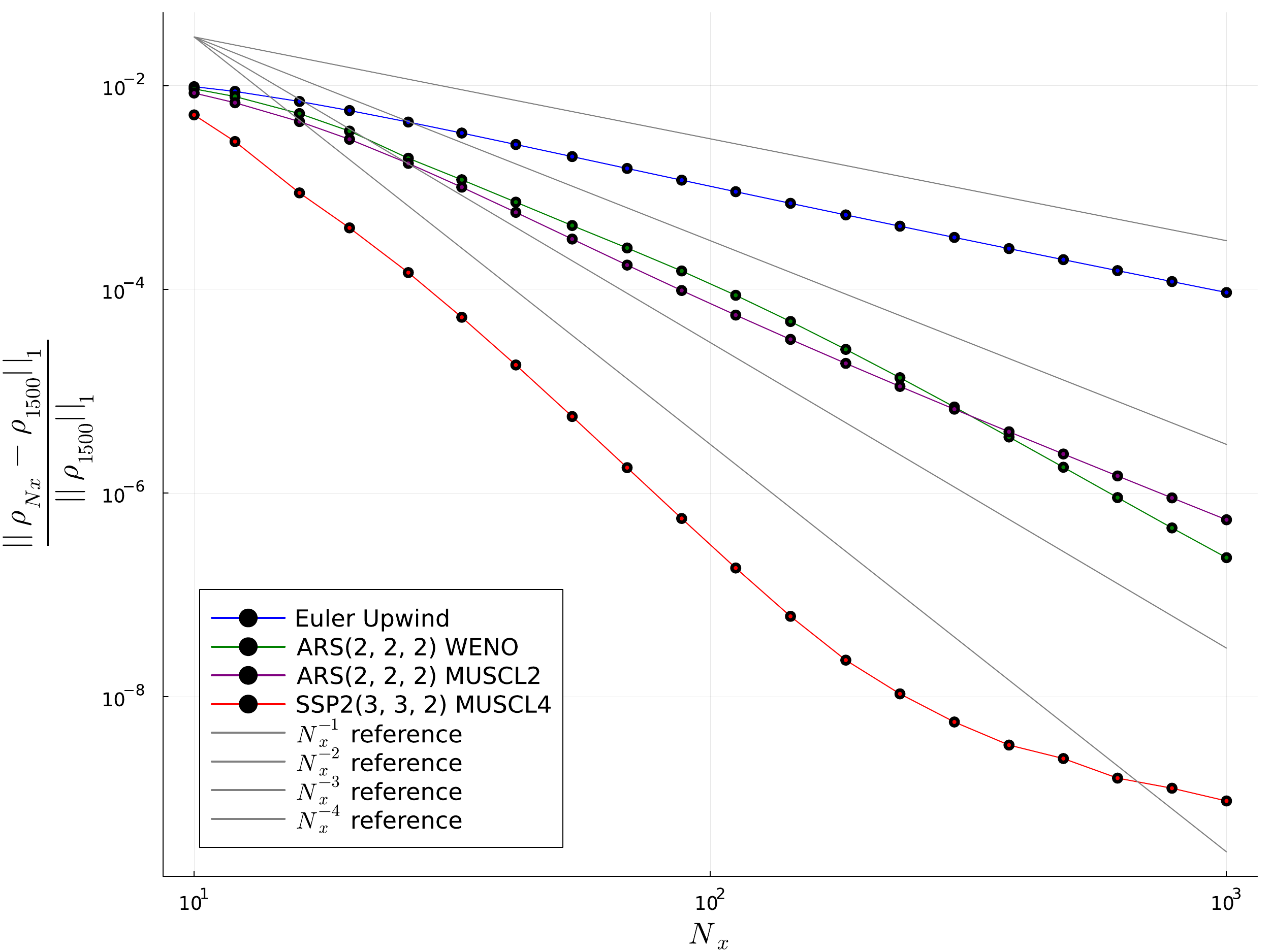}
		\label{fig:1DConvergenceRho}
	\end{subfigure}
	\begin{subfigure}[b]{0.49\textwidth}
		\centering
		\includegraphics[width=\linewidth]{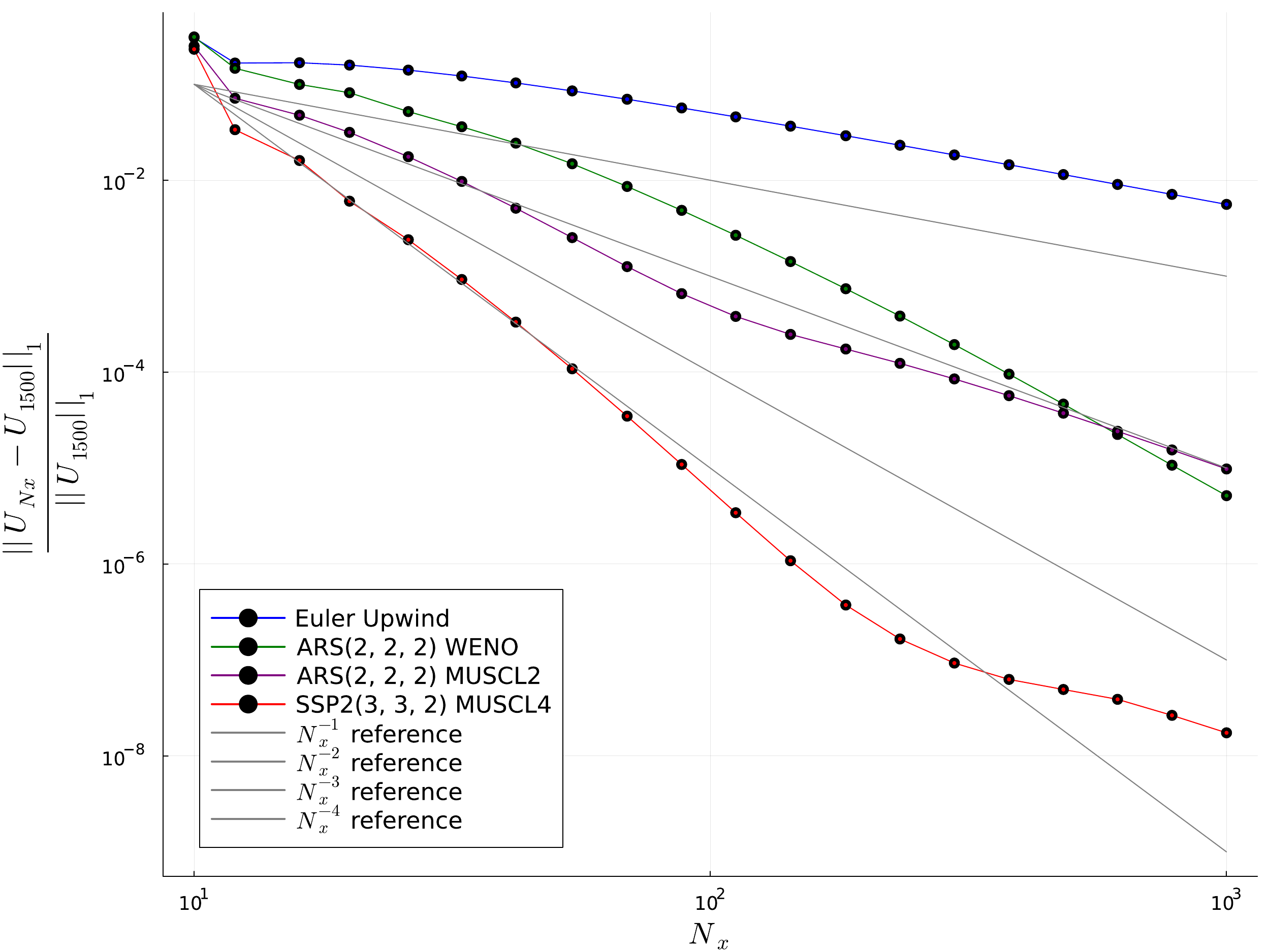}
		\label{fig:1DConvergenceU}
	\end{subfigure}
	\caption{The relative \(L^1\) error of the density (left) and mean velocity (right) for several schemes for a smooth initial condition in one spatial dimension.}
	\label{fig:1DConvergence}
\end{figure}

The results are plotted in figure \ref{fig:1DConvergence}. The first-order positive fallback scheme from section \ref{section:Space1DFirstOrder} is labelled `Euler Upwind'. It achieves first order as advertised. The WENO MLS method achieves third-order accuracy in this test case, although it is only second-order. The second-order MUSCL method achieves similar performance as the WENO MLS method. The fourth-order MUSCL method achieves fourth order accuracy, after which the error degrades to second-order. This is the time integration error of the second-order SSP2(3, 3, 2) method that starts to dominate the total error. 

\begin{remark}
	Since we fix the number of grid points in velocity space, we are reporting on the convergence of a discrete velocity model, since there is some error related to the velocity discretisation that we are not measuring here. For small Knudsen numbers, as is the case here, typically only a small amount of grid points is necessary, and the error due to the space discretisation is usually more relevant. 
\end{remark}

\subsection{Convergence in two space dimension}
\begin{figure}[h!]
	\centering
	\begin{subfigure}[b]{0.49\textwidth}
		\centering
		\includegraphics[width=\linewidth]{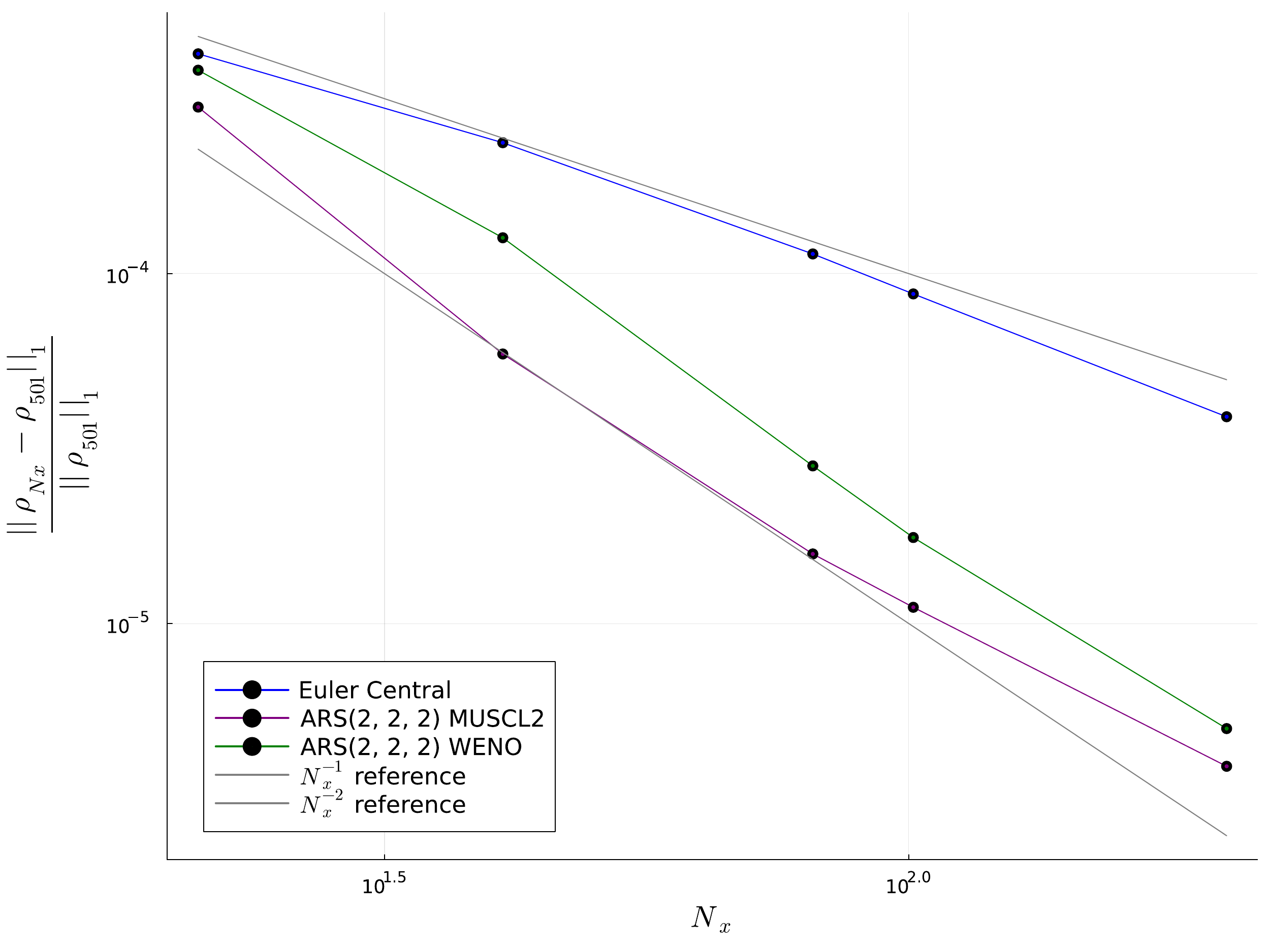}
		\label{fig:2DConvergenceRho}
	\end{subfigure}
	\begin{subfigure}[b]{0.49\textwidth}
		\centering
		\includegraphics[width=\linewidth]{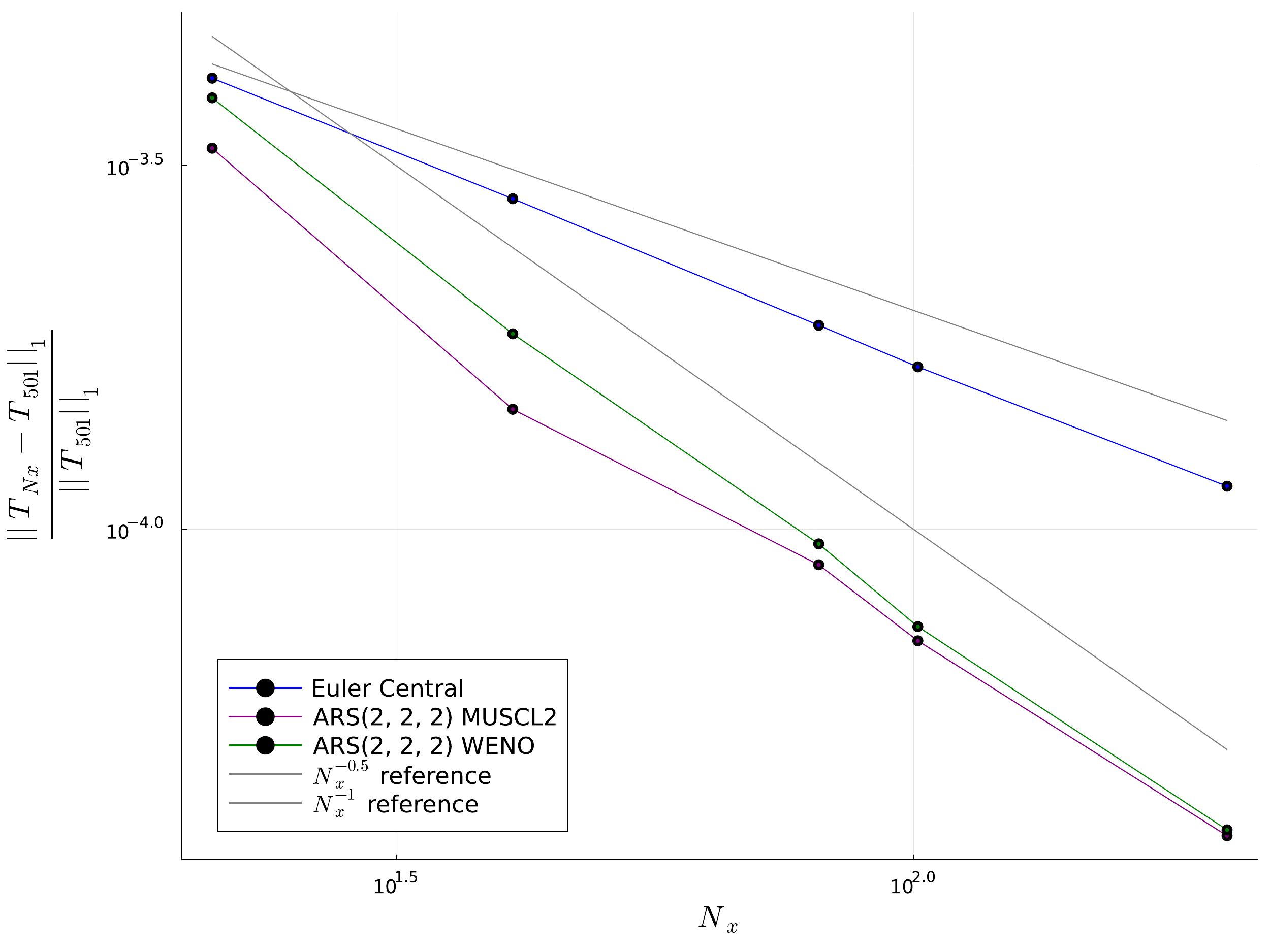}
		\label{fig:2DConvergenceT}
	\end{subfigure}
	\caption{The relative \(L_1\) error of the density (left) and temperature (right) for several schemes for a smooth initial condition in two spatial dimension. Because the convergence of the mean velocity is the same as the convergence for the density for all methods, the plots are omitted.}
	\label{fig:2DConvergence}
\end{figure}

Similarly to the one dimensional case, we perform a convergence test in two space dimensions. We create a fixed \([-1, 1]^2\) domain with diffuse reflective boundary conditions with wall temperature \(T_w = 1\). Inside the domain is a gas with specific gas constant \(R_s = 1\) and relaxation time \(\tau = 10^{-5}\). Initially, the gas is in local thermodynamic equilibrium with the density and temperature set to one and the mean velocity set to \begin{align}
	U_x^0 &= \frac{1}{10}\left[ \exp \left\{ -\left( 10 \sqrt{ \left( x - 0.2 \right)^2 + y^2 } - 1 \right)^2 \right\} - 2 \exp \left\{- \left( 10 \sqrt{ \left( x + 0.2 \right)^2 + y^2 } - 1 \right)^2 \right\}  \right]\\
	U_y^0 &= \frac{1}{10}\left[ \exp \left\{ -\left( 10 \sqrt{ \left( y - 0.2 \right)^2 + x^2 } - 1 \right)^2 \right\} - 2 \exp \left\{- \left( 10 \sqrt{ \left( y + 0.2 \right)^2 + x^2 } - 1 \right)^2 \right\}  \right].
\end{align} The velocity variable is discretised using \(N_v = 30^2\) uniformly distributed points in the domain \([-v_{\rm max}, v_{\rm max}]^2\) with \(v_{\rm max} = 15\). We perform the simulation up to time \(t_{f} = 5 \times 10^{-3}\) with CFL number \(0.02\). We set the time step \(\Delta t\) according to a ``classical'' fixed uniform grid time step restriction 
\begin{align}
\Delta t = \frac{2 \text{CFL}}{N_x v_{\rm max}}, \label{eq:2DConvergenceTimeStep}
\end{align} 
This is sufficient to ensure stability throughout the simulation. We compute the \(L^1\) error of the numerical solution with respect to a fine simulation computed with \(501^2\) grid points. The errors are computed by interpolating the numerical solution to a uniformly distributed fixed grid of \(501^2\) points using a cubic interpolation routine. The convergence plot of the density and temperature are given in figure \ref{fig:2DConvergence}. The first-order positive scheme from section \ref{section:Space2DFirstOrder} is labelled ``Euler Central''. For comparison purposes, we have also implemented the second-order meshless WENO method from \cite{tiwari_meshfree_2022}, that is based on the WENO-SPH method from \cite{avesani_new_2014}. All methods achieve the expected order for the density. Because the solution is smooth, the meshless MUSCL method with MOOD yields the same result as the same method without MOOD. For the temperature, the order of all methods is decreased.  

\subsection{Shock tube}
\label{section:shockTube}
In this subsection, and the next two, we shall make use of dimensional quantities and adopt the SI system, so that lengths are expressed in meters, times in seconds, masses in Kilograms, pressures in Pascals, temperatures in degrees Kelvin, and so on. The monoatomic gas we consider is Argon, and has a specific gas constant equal to \(R_s = 208 \) J/Kg K.

\begin{figure}
	\centering
	\begin{subfigure}[b]{0.49\textwidth}
		\centering
		\includegraphics[width=\linewidth]{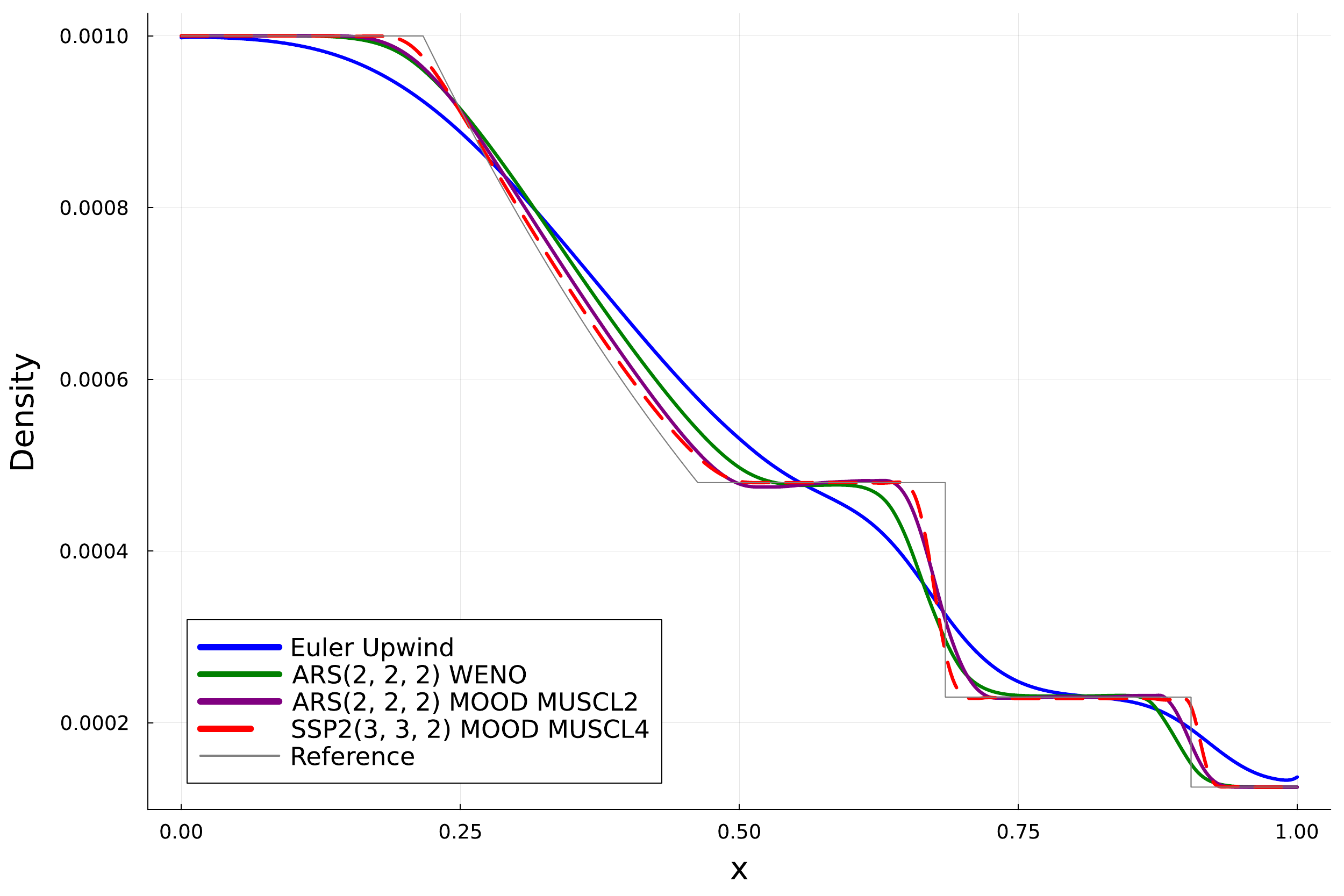}
		\label{fig:riemannRho}
	\end{subfigure}
	\begin{subfigure}[b]{0.49\textwidth}
		\centering
		\includegraphics[width=\linewidth]{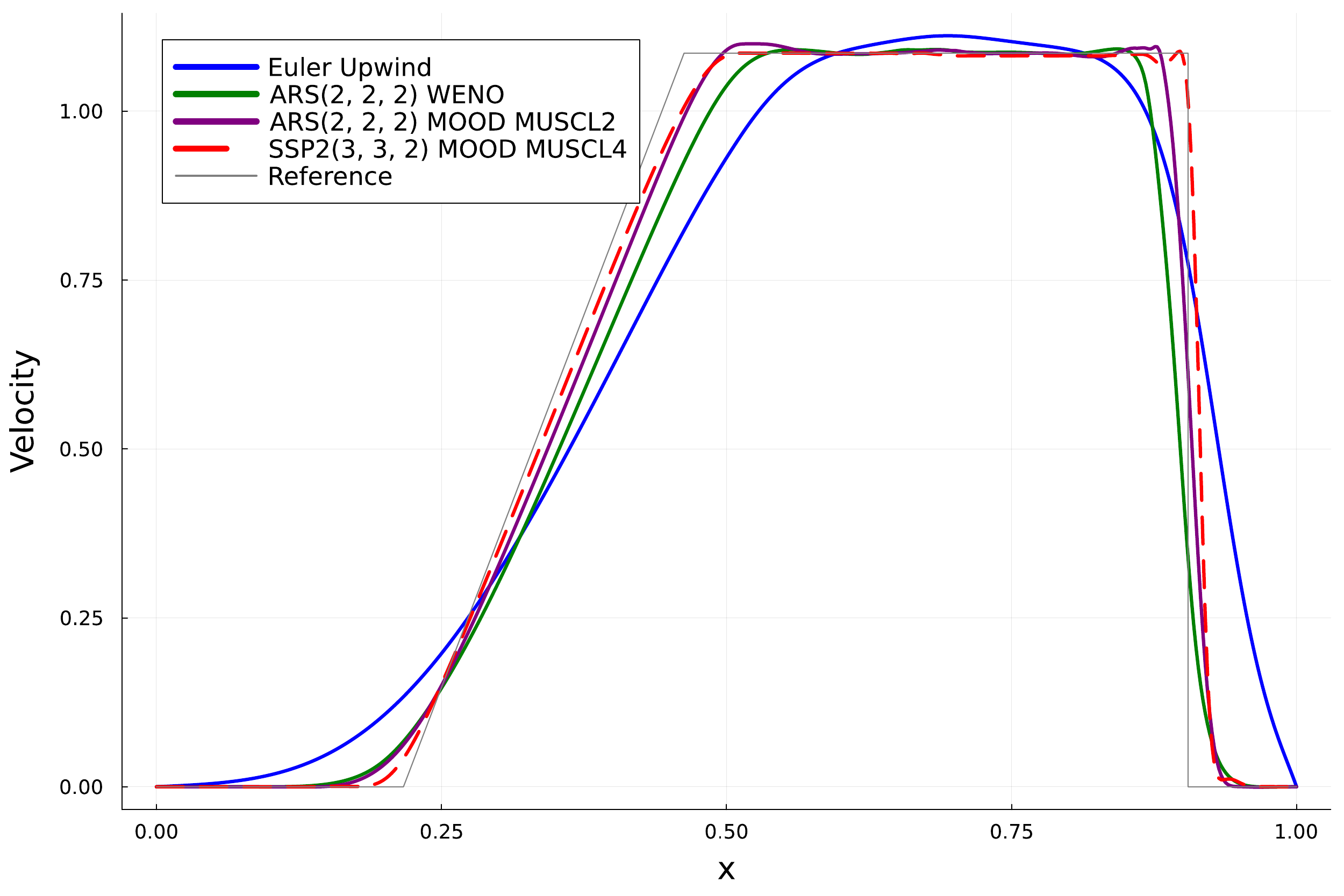}
		\label{fig:riemannU}
	\end{subfigure}
	\begin{subfigure}[b]{0.49\textwidth}
		\centering
		\includegraphics[width=\linewidth]{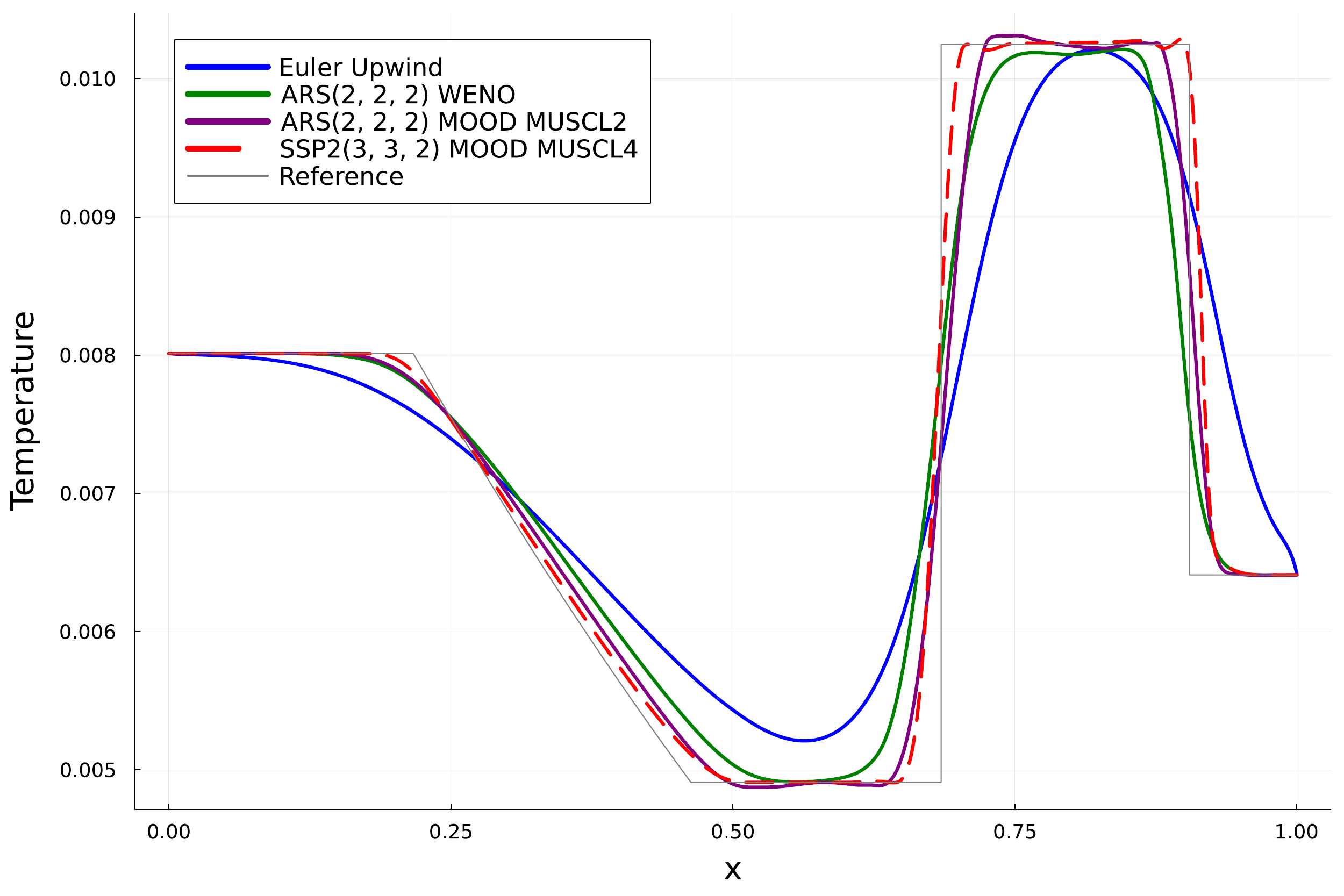}
		\label{fig:riemannT}
	\end{subfigure}
	\caption{The numerical solution and the exact solution to Sod's shock tube at time \(0.17 \;\SI{}{s}\). The numerical solution is computed using \(N_x = 100\), \(N_v = 80\) and CFL = \(0.5\).}
	\label{fig:1DShockTube}
\end{figure}

To illustrate the correct behaviour of MOOD, we simulate Sod's shock tube. We consider a fixed spatial domain \([0, 1]\) that is discretised using 100 points. The velocity variable is discretised using 80 points in the interval \([-20, 20]\). Initially, the gas is in thermodynamic equilibrium with macroscopic quantities \begin{align}
	(\rho, U, T) = \begin{cases}
		\rho_l = 1\times10^{-3} \; \SI{}{kg/m^3}, \; U_l = 0 \;\SI{}{m/s}, \;T_l = 8.012 \times 10^{-3} \; \SI{}{K} & \text{if \(x \leq 0.5 \; \SI{}{m}\)} \\
		\rho_r = \frac{1\times10^{-3}}{8}\; \SI{}{kg/m^3}, \;U_r = 0 \;\SI{}{m/s}, \;T_r = 6.41 \times 10^{-3} \; \SI{}{K} & \text{else. }
	\end{cases}
\end{align} We set the relaxation time \(\tau\) to \(10^{-6} \; \SI{}{s}\) such that we can compare the numerical solution to the exact solution of the Riemann problem of the compressible Euler equations. The simulation is performed with several algorithms with CFL number 0.5 up to time \(t_f = 0.17 \; \SI{}{s}\). The results are plotted in figure \ref{fig:1DShockTube}. The first-order fallback scheme yields a very diffusive result. The result using the WENO MLS scheme is somewhat diffusive, but does not contain any oscillations. The MUSCL methods yield the sharpest result at shocks, however, the fourth order scheme produces some minor oscillations. These oscillations were also observed for the linear advection equation \cite{willems2025b} with shock initial condition.

\begin{figure}
    \centering
    \begin{subfigure}[b]{0.49\textwidth}
	\centering
	\includegraphics[width=0.75\linewidth]{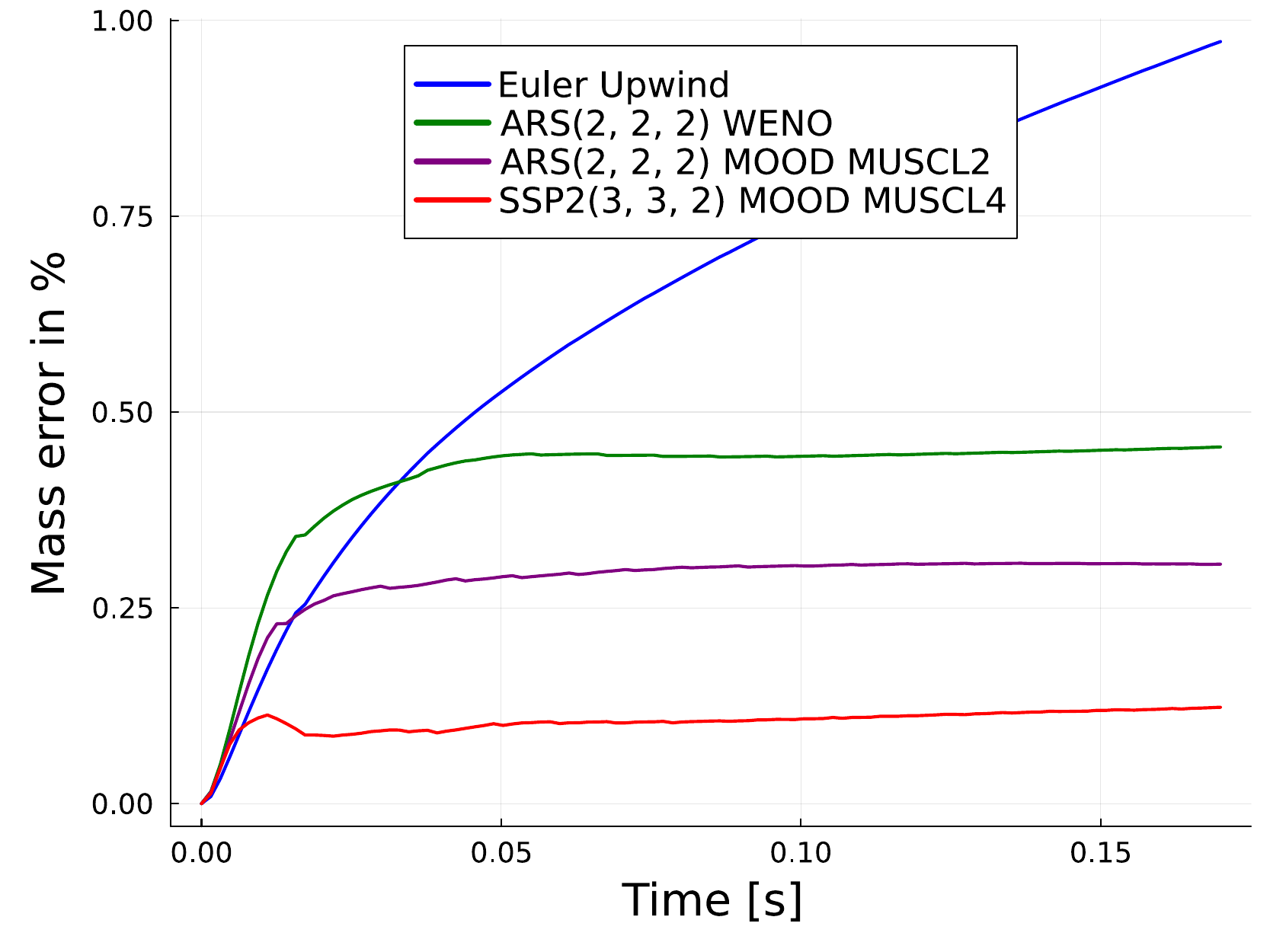}
	\subcaption{The mass error throughout time expressed as a percentage.}
	\label{fig:massLoss}
    \end{subfigure}
    \begin{subfigure}[b]{0.49\textwidth}
	\centering
	\includegraphics[width=0.75\linewidth]{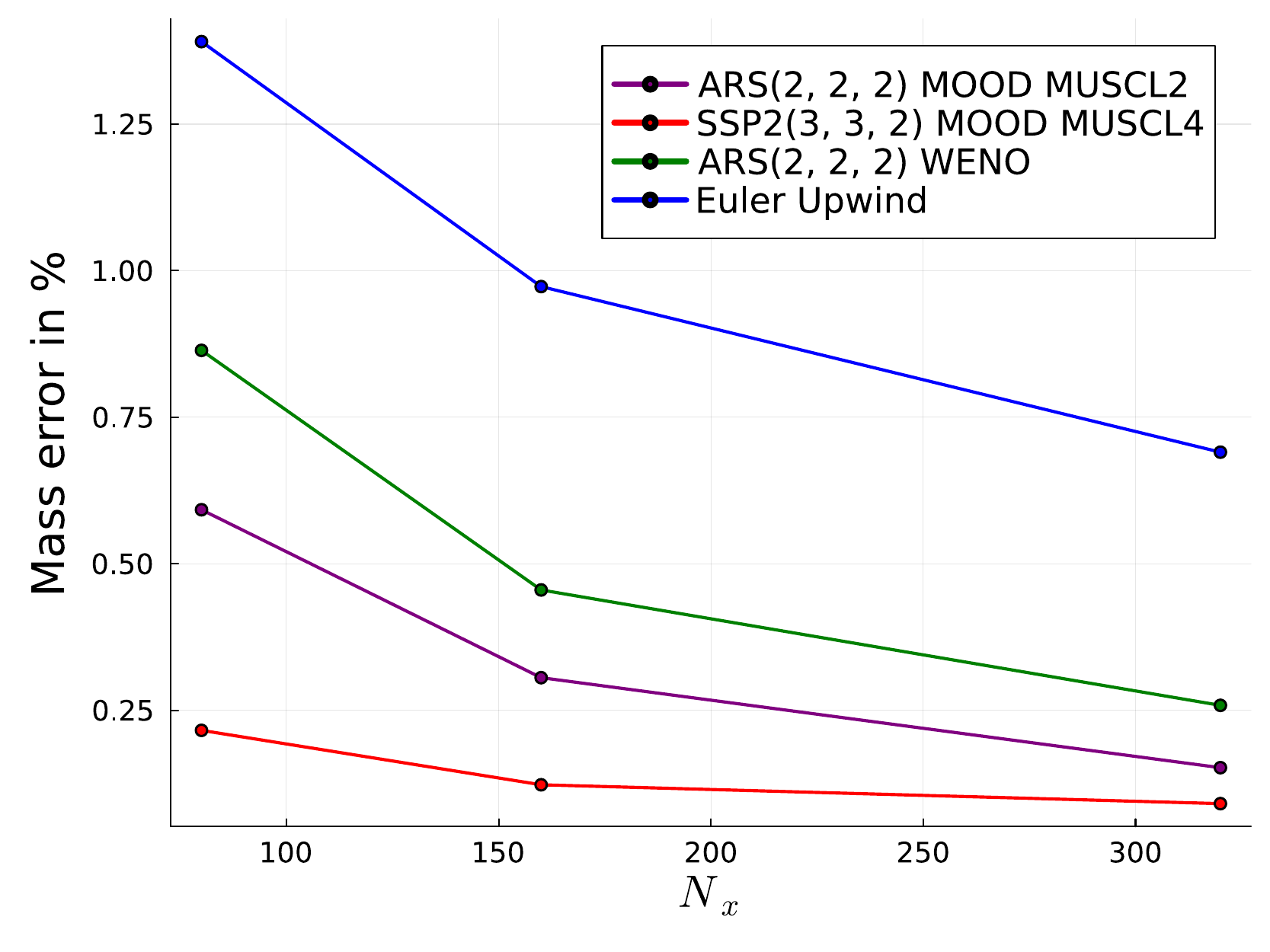}
	\subcaption{The mass error at final time \(0.17 \; \SI{}{s}\) as a function of the amount of grid points.}
	\label{fig:massLossConvergence}
    \end{subfigure}
    \caption{The mass error for several numerical methods for Sod's shock tube. }
\end{figure}

To show the lack of conservation of the meshless methods, we repeat the simulation with \(N_x=80, 160\) and \(320\). All other simulation parameters are fixed as in the previous simulation. The error on the mass is computed throughout the simulation and expressed as a percentage with respect to the mass of the initial condition. The error in the mass for \(N_x=80\) is plotted in figure \ref{fig:massLoss}. The first-order fallback method from \ref{section:Space1DFirstOrder} yields an error that increases rather quickly throughout time up to \(\approx 1.0\%\) at the end of the simulation. The WENO and MUSCL methods initially add mass to the system and then stabilize. The error of the mass is clearly related to the accuracy of the numerical method. In figure \ref{fig:massLossConvergence}, we plot the error at the final time for different grid sizes. We observe that the error decreases as the grid is refined. 

\subsection{The moving plate problem}
\begin{figure}
	\centering
	\includegraphics[width=0.7\linewidth]{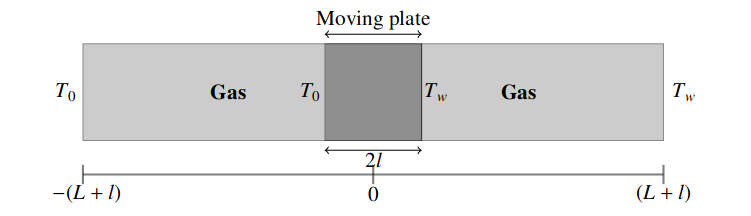}
	\caption{Schematic of the moving plate problem with \(L = 1 \; \SI{}{m}\), \(l = 0.1 \; \SI{}{m}\), \(T_0 = 270 \; \SI{}{K}\) and \(T_w = 330 \; \SI{}{K}\).}
	\label{fig:slider}
    
\end{figure}

We consider the moving plate problem from \cite{dechriste_numerical_2012}. Two gaseous domains are separated by a moving plate of size \(l=0.1 \; \SI{}{m}\). Initially, the gas (\(R_s = 208 \>{\rm J}/({\rm Kg~}^\circ {\rm K})\)) in both chambers is stationary and in thermodynamic equilibrium with pressure \(p_0 = 3.86\times10^{-2} \; \SI{}{Pa}\) and temperature \(T_0=270 \; \SI{}{K}\). The relaxation time of the gas is computed from \cite{tiwari_meshfree_2022, chapman_mathematical_1999} \begin{align}
	\lambda = \frac{k_B}{\sqrt{2} \pi \rho_0 R_s d^2} \approx 1.605 \times 10^{-1} \; \SI{}{m}, \quad \tau = \frac{4 \lambda}{\sqrt{8 \pi R_s T_0}} \approx 5.404\times10^{-4} \; \SI{}{s} \label{eq:relaxationTime}, 
\end{align} 
and assumed constant throughout simulation, where \(k_B\) is the Boltzmann constant, \(\lambda\) the mean free path, \(\rho_0\) the initial density of the gas and \(d = 3.68\times10^{-10} \; \SI{}{m}\) the diameter of a gas molecule. The density of the moving plate is chosen \(\rho_p = 10 \rho_0\). The walls of the left and right chamber are heated to \(270 \; \SI{}{K}\) and \(330\; \SI{}{K}\) respectively. This causes the pressure of the gas in the right chamber to heat, which in turn, pushes the moving plate to the left. The force the gas applies on the moving slider in this simplified one dimensional settings is \begin{align}
	F = (\phi_l - \phi_r)A, \label{eq:EqSlider}
\end{align} with \(\phi_l\) and \(\phi_r\) the pressure of the gas at the left and right boundary, and \(A = 2l\rho_p\) the surface area of the moving plate. Due to the one dimensional nature of the simulation, the plate does not rotate. The equilibrium state of the moving plate is \begin{align}
	x_{eq} = L \frac{T_0 - T_w}{T_0 + T_w} = -0.1 \; \SI{}{m}.
\end{align}
The gaseous domain is discretised using 54 points, two points of which lie exactly on the boundary of the moving plate and move with the velocity of the plate. The velocity variable is discretised using 20 uniformly distributed points in the domain \([-1400, 1400]\). The simulation is performed using several meshless methods with the CFL number set to 0.5 up to time \(t_f = 0.2 \; \SI{}{s}\). The equation of motion of the slider \eqref{eq:EqSlider} is discretised using an explicit Euler time integration method. Although this could potentially lead to a loss of accuracy, based on numerical testing, the discretisation error of the gas simulation is the dominant one. 

\begin{figure}[h!]
	\centering
	\begin{subfigure}[b]{0.49\textwidth}
		\centering
		\includegraphics[width=\linewidth]{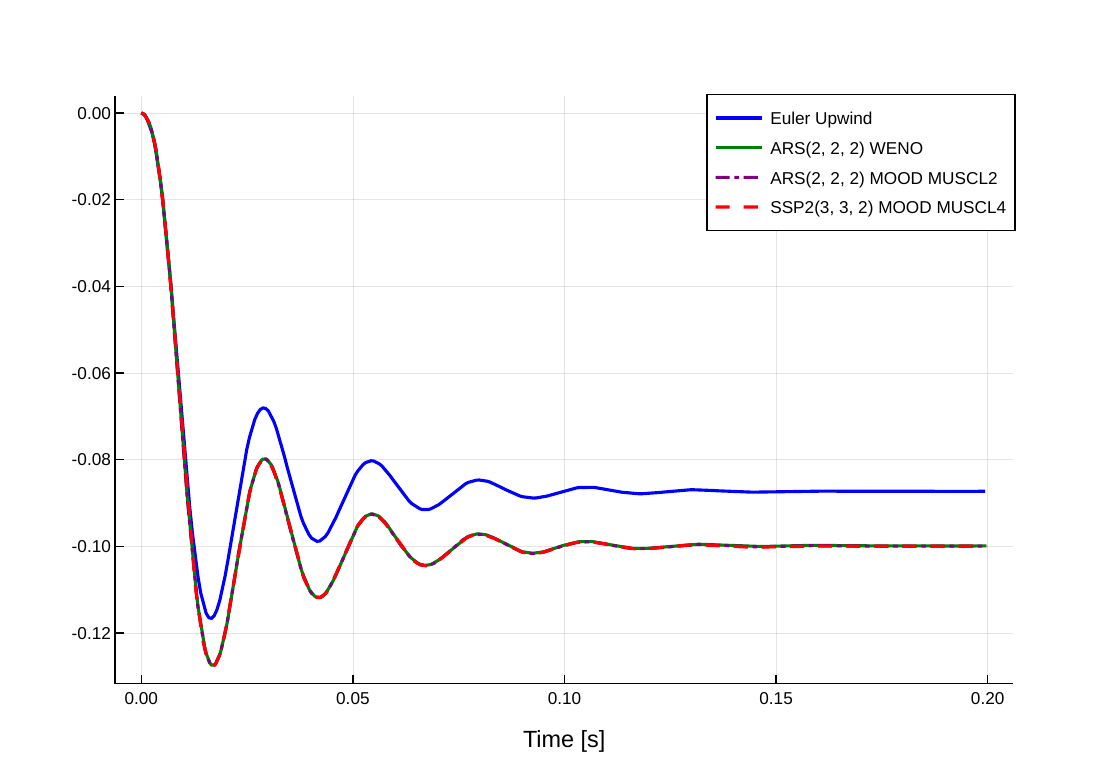}
		\label{fig:sliderPosition}
	\end{subfigure}
	\begin{subfigure}[b]{0.49\textwidth}
		\centering
		\includegraphics[width=\linewidth]{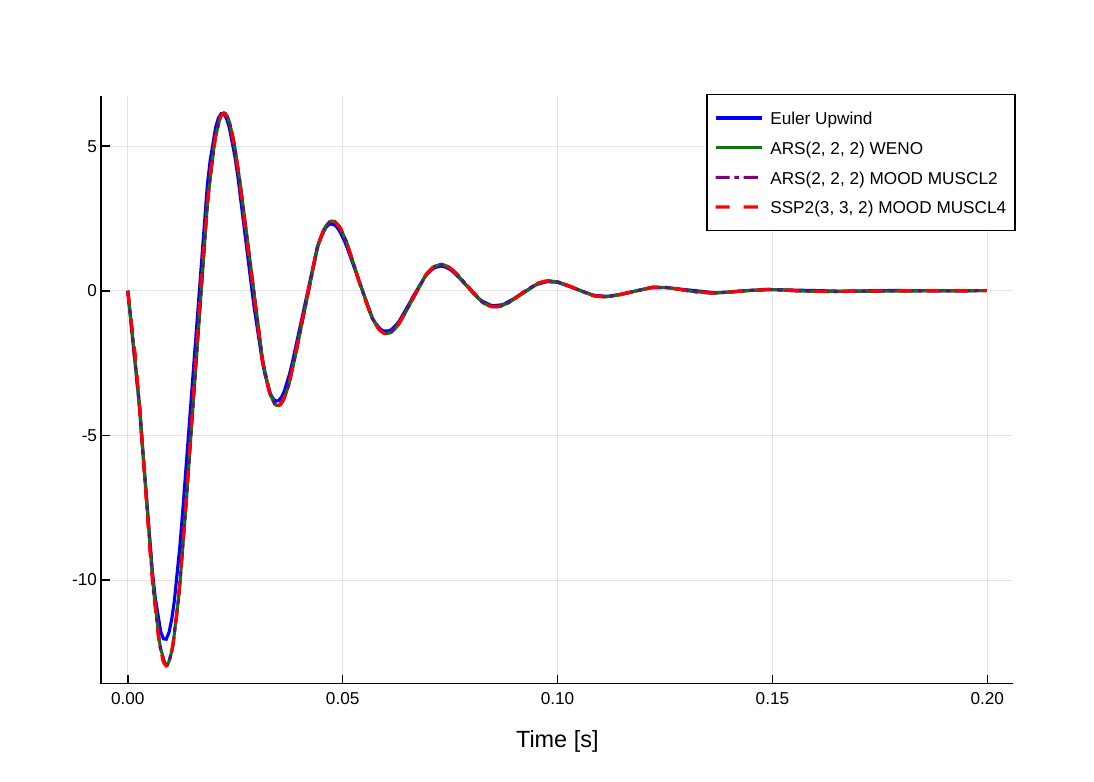}
		\label{fig:sliderVelocity}
	\end{subfigure}
	\caption{The position (left) and velocity (right) of the moving plate up to time \(t_f = 0.2\).}
	\label{fig:sliderNumerics}
\end{figure}

The time step for integrating the position of the moving plate is chosen to match the time index of the Runge-Kutta stages of the scheme used to advance the points within the domain. In this way, points on the boundary and inside the domain are always at the same time index. The position and velocity of the moving plate throughout time is plotted in figure \ref{fig:sliderNumerics}. The first-order method yields the wrong equilibrium solution of the slides, while the higher-order methods all perform equally well. Note that due to the moving nature of the grid, the grid points automatically adapt to the shape of the domain. The method is thus immediately applicable to simulations with moving boundaries with the need for a cut-cell method \cite{dechriste_cartesian_2016} or immersed boundary method \cite{dechriste_numerical_2012}.

\subsection{The driven square cavity}
In this test, we demonstrate the ability of the meshless method to simulate arbitrary rigid bodies immersed in the rarefied gas using the driven square cavity test case \cite{tiwari_meshfree_2022, banda_lattice-boltzmann_2007}. The driven square cavity problem consists of a closed box \([0, L]^2, L = 1\times10^{-6} \; \SI{}{m}\). The top boundary of the box moves horizontally with velocity \begin{align}
	U_\text{wall}^x = U_\text{max} \left(\frac{2x}{L}\right)^2 \left( 2 - 2 \frac{x}{L} \right)^2,
\end{align} with \(U_\text{max} = 10 \; \SI{}{m/s}\). We apply diffuse-reflective boundary conditions with temperature \(T_w = 270 \; \SI{}{K}\) on all four walls of the box. The box contains argon gas, which is initially at rest and in thermodynamic equilibrium at a a temperature of \(T^0 = 270 \; \SI{}{K}\). 

We perform the simulation with two initial pressures, which in turn leads to different Knudsen numbers; see table \ref{tab:squareCavityParam}. For the first pressure, the gas is in the slip flow regime; for the second pressure, the gas is in the transition regime. 
\begin{table}[h!]
\centering
\begin{tabular}{|c||c|c|c|c|c|c|} 
	\hline
	& \(p\) & \(\tau\) & \(\lambda\) & \({\rm Kn}\) & \({\rm Re}\) & \({\rm Ma}\) \\ \hline \hline
	1 & \(61776 \; \SI{}{Pa}\) & \(3.37671 \times 10^{-11} \; \SI{}{s}\) & \(1.00292 \times 10^{-8} \; \SI{}{m}\) & \(0.01\) & \(5.3\) & \(0.036\) \\
	2 & \(617760 \; \SI{}{Pa}\) & \(3.37671 \times 10^{-10} \; \SI{}{s}\) & \(1.00292 \times 10^{-7} \; \SI{}{m}\) & \(0.1\) & \(0.53\) & \(0.036\) \\
	\hline
\end{tabular}
\caption{From left to right, \(p^0 = \rho^0 R_s T^0\) is the initial pressure, \(\tau\) is the mean free time, \(\lambda\) is the mean free path of which the definition is given in \cite{tiwari_meshfree_2022}, \(Kn\) the Knudsen number and \({\rm Re}\) the Reynolds number. 
The latter two are defined as  \({\rm Kn} = {\lambda}/{L}\) and \({\rm Re} = {\rho^0 U_{\rm max} L}/{\tau p^0}\).}
\label{tab:squareCavityParam}
\end{table} First, we perform a simulation of the gas in slip flow regime and compare it with several other methods. We consider a first-order and second-order semi-Lagrangian method \cite{groppi_high_2014}, Direct Simulation Monte Carlo (DSMC) \cite{bird_molecular_1994}, the first-order meshless method presented in section \ref{section:Space2DFirstOrder}, and the second-order meshless MUSCL method presented in section \ref{section:Space2DSecondOrder}. We note that with DSMC, we refer to a Monte Carlo particle method associated to the Boltzmann equation. Although, this is a different model to describe rarefied gases, in the regimes considered here, we do not observe any differences.

\begin{figure}[h!]
	\centering
	\includegraphics[width=0.75\linewidth]{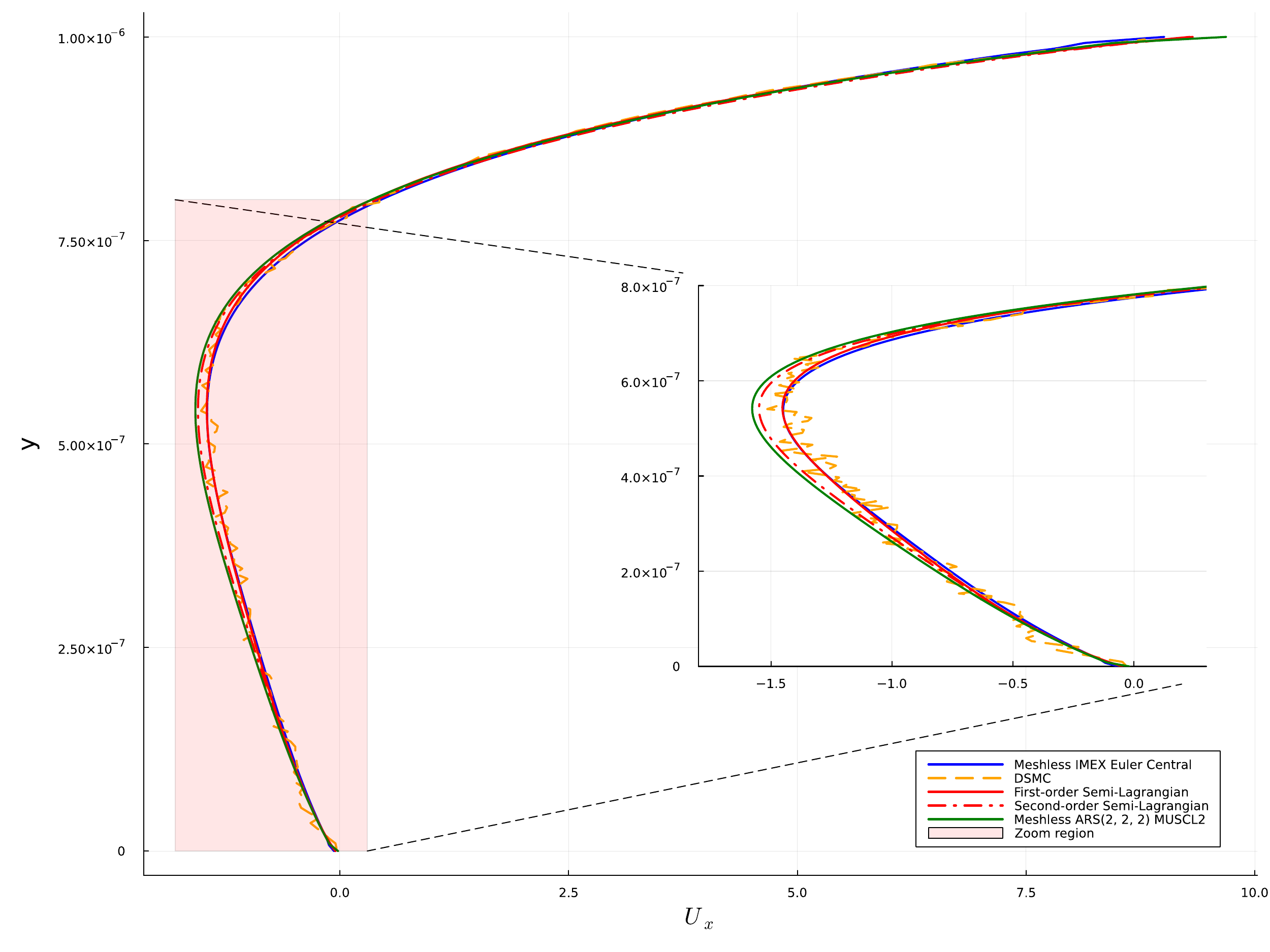}
	\caption{The velocity in the $x$-direction along the vertical centre line of the square cavity. The gas is in slip flow regime, see the first row in table \ref{tab:squareCavityParam}.}
	\label{fig:SquareCavityCentreLineComparison}
\end{figure}

For the deterministic methods, we use a fine grid of \(120^2\) grid points (initially uniform), a time step of \(\Delta t = 1 \times 10^{-12}\), and perform a simulation up to time \(t_f = 8 \times 10^{-8} \; \SI{}{s}\). Specifically for the BGK equation, we use a velocity grid of \(N_v^2 = 30^2\) with \(v_{\rm max} = 1500\). For the DSMC simulation, we initially generated 50 particles per grid square (same resolution as above), and averaged over approximately 2 billion time steps after steady state is reached. The velocity in the $x$-direction is plotted in figure \ref{fig:SquareCavityCentreLineComparison} along the centre line of the square cavity.
 We observe good agreement between the results obtained from the first-order methods (DSMC, first-order semi-Lagrangian, first-order meshfree method), and the second-order methods (second-order semi-Lagrangian and meshfree method). We conclude that the new implementation of the diffuse-reflective boundary conditions is correct. \begin{comment}
	We note that the methods presented here for the BGK equation are not asymptotic preserving in the low-Mach limit. 
\end{comment}

\begin{figure}[h!]
	\centering
	\begin{subfigure}[b]{0.49\textwidth}
		\centering
		\includegraphics[width=\linewidth]{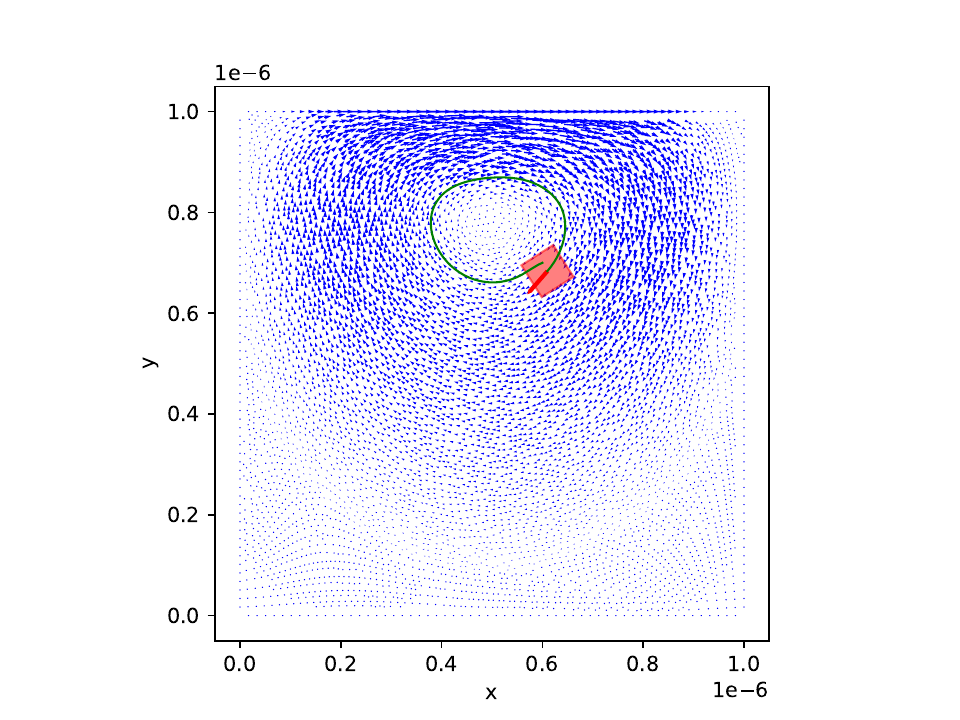}
		\label{fig:sim1quiver500}
		\subcaption{\(Kn = 0.01, t = 1\times10^{-6}\) }
	\end{subfigure}
	\begin{subfigure}[b]{0.49\textwidth}
		\centering
		\includegraphics[width=\linewidth]{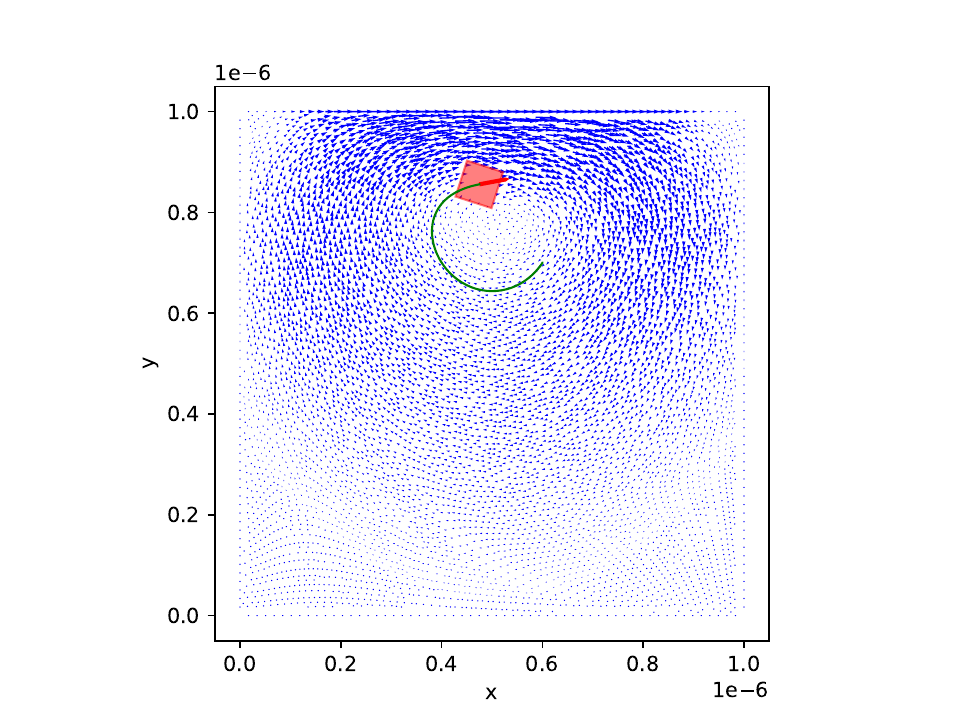}
		\label{fig:sim2quiver500}
		\subcaption{\(Kn = 0.1, t = 5\times10^{-7}\) }
	\end{subfigure}
		\begin{subfigure}[b]{0.49\textwidth}
		\centering
		\includegraphics[width=\linewidth]{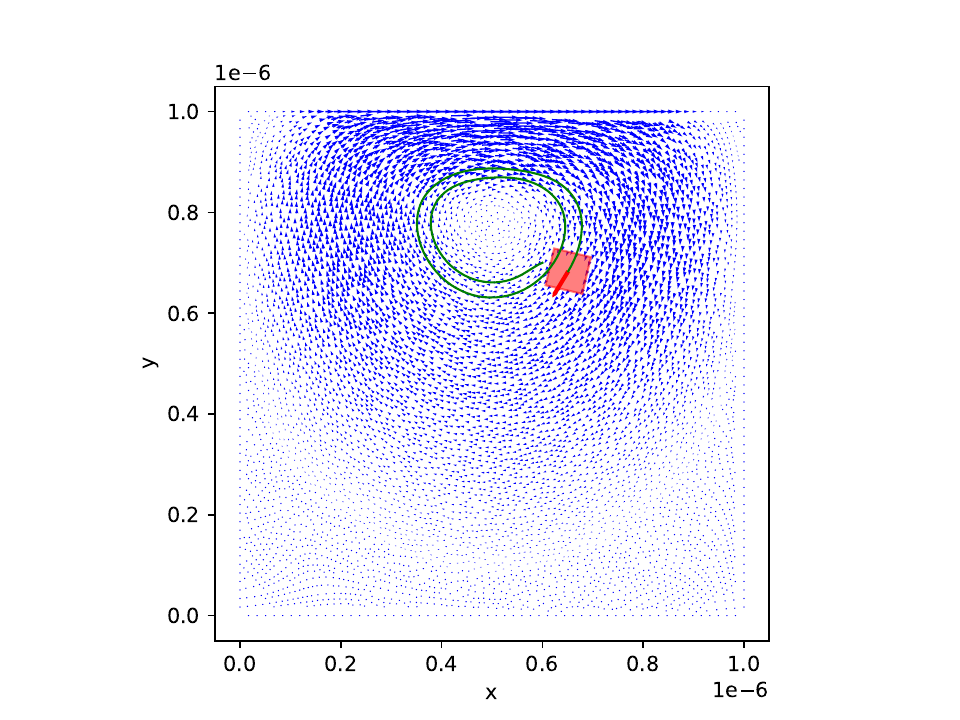}
		\label{fig:sim1quiver1000}
		\subcaption{\(Kn = 0.01, t = 5\times10^{-7}\) }
	\end{subfigure}
	\begin{subfigure}[b]{0.49\textwidth}
		\centering
		\includegraphics[width=\linewidth]{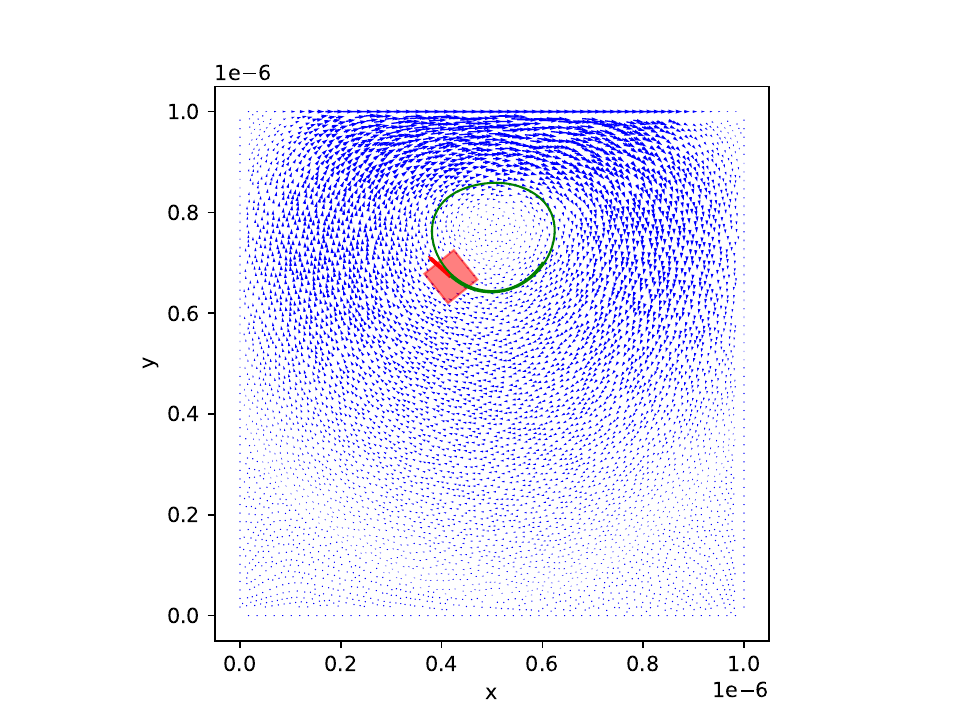}
		\label{fig:sim2quiver1000}
		\subcaption{\(Kn = 0.1, t = 1\times10^{-6}\) }
	\end{subfigure}
	\caption{Velocity field of the driven square cavity simulation plotted at several times for various Knudsen numbers. The square rigid body is plotted in red. The path of the centre of mass of the rigid body is plotted in green. }
	\label{fig:squareCavity}
\end{figure}

Next, at position \((6\times10^{-7}, 7\times 10^{-7})\), we place a square rigid body with density \(\rho_{rb} = 10\rho^0\) and temperature \(T_{rb} = 270 \; \SI{}{K}\). We now perform two simulations, of which the initial pressures are given in table \ref{tab:squareCavityParam}. We perform the simulations using \(60^2\) spatial discretisation points. We discretise the velocity variable using \(N_v = 30\) points in both velocity dimensions and set \(v_{\rm max} = 1500\). We fix 16 additional grid points to the border of the rigid body to accurately capture its motion. We perform the simulation up to time \(t_f = 1\times10^{-6} \; \SI{}{s}\) with time step \(\Delta t = 5\times10^{-12} \). The results are plotted in figure \ref{fig:squareCavity}. 

\subsection{Shear layer problem}
As a final example, we consider the shear layer problem from \cite{banda_lattice-boltzmann_2007, kurganov_third-order_2000, bell_second-order_1989}. Note that in this example, we do not consider physical parameters. We use a periodic domain \([0, 2\pi]^2\), initially discretised with \(80^2\) uniformly distributed points. Inside the domain is a gas with gas constant \(R_s = 1\). Initially, the gas is in thermodynamic equilibrium, with density \(\rho^0 = \frac{15}{\pi}\), mean velocity \begin{align}
	U_x^0 &= \begin{cases}
		\tanh \left(\frac{15}{\pi} \left(y -  \frac{\pi}{2}\right) \right), \quad y \leq \pi \\
		\tanh \left(\frac{15}{\pi} \left(\frac{3\pi}{2} - y \right) \right), \quad y > \pi
	\end{cases} \\
	U_y^0 &= 0.05 \sin(x), 
\end{align} and temperature \(T^0 = \frac{1}{2\rho^0}\). The velocity variable is discretised using \(N_v^2 = 25^2\) points in the interval \([-6\sqrt{R_s T^0}, 6\sqrt{R_s T^0}]^2\). The mean free time, time step and several non-dimensional parameters are given in table \ref{tab:shearFlowParam}.

\begin{table}[h!]
	\centering
	\begin{tabular}{|c|c|c|c|c|} 
		\hline
		\(\tau\) & \(\Delta t\) & \({\rm Kn}\) & \({\rm Ma}\) & \(Re\) \\ \hline
		\(1 \times 10^{-5}\) & \(1.17617 \times 10^{-3}\) & \(1 \times 10^{-5}\) & \(0.1585\) & \(26318.9\) \\
		\hline
	\end{tabular}
	\caption{The mean free time \(\tau\), time step \(\Delta t\), the Knudsen number \({\rm Kn}\), the Mach number \({\rm Ma}\), and the Reynolds number \(Re\). The chosen time step yields a CFL number of 0.45 for the semi-Lagrangian method.}
		
	\label{tab:shearFlowParam}
\end{table}

The simulation is performed using the second-order meshless MUSCL method with a the ARS(2, 2, 2) time integration method without MOOD. We compare the numerical solution of the meshless method with a second-order in time, third-order in space semi-Lagrangian scheme on a Cartesian mesh. The reconstruction in space is performed using a bicubic interpolation procedure. 

\begin{figure}[h!]
	\centering
	\begin{subfigure}[b]{0.49\textwidth}
		\centering
		\includegraphics[width=\linewidth]{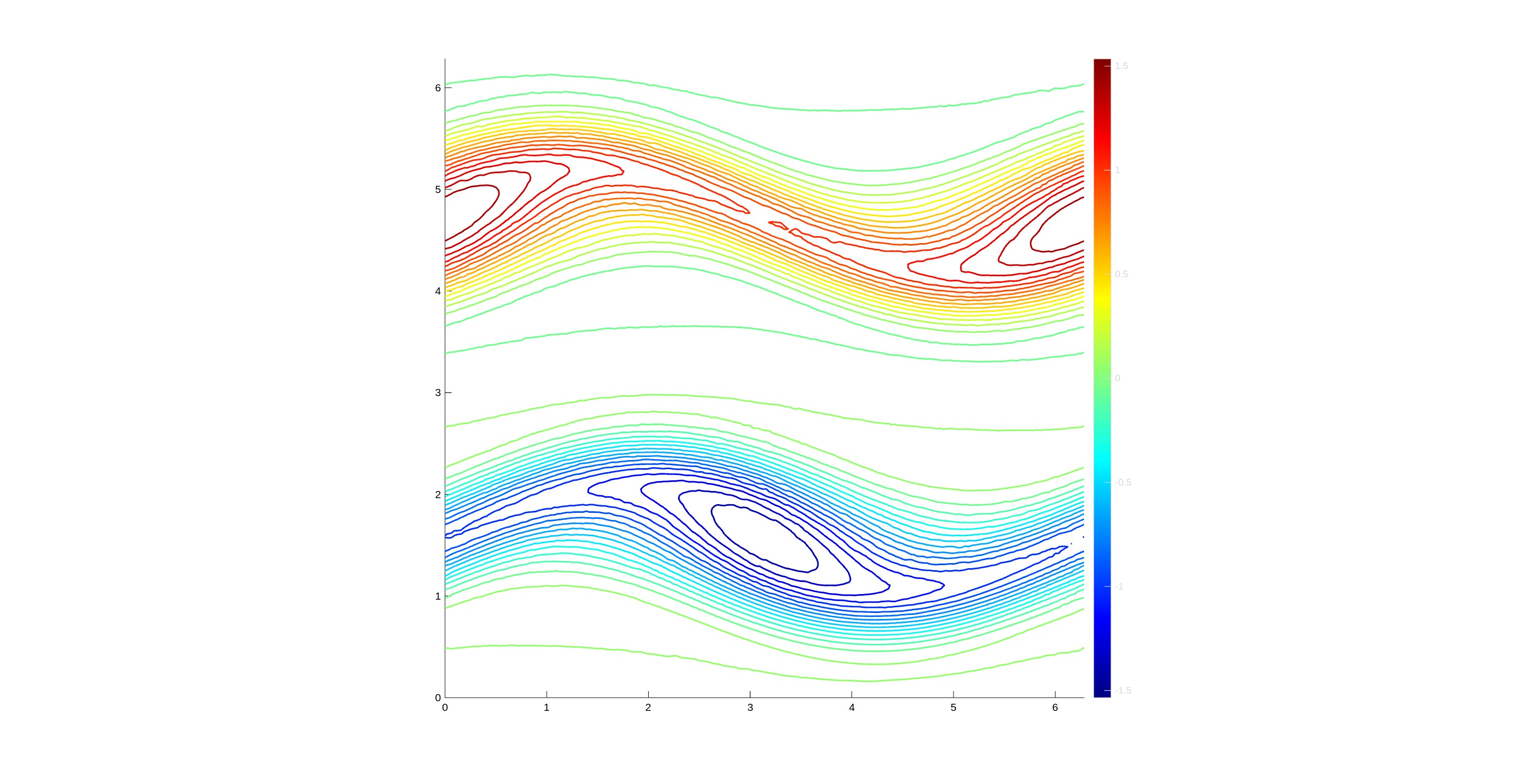}
		\caption{Meshfree MUSCL method, time 4.93993.}
		\label{fig:MUSCLVorticity5}
	\end{subfigure}
	\begin{subfigure}[b]{0.49\textwidth}
		\centering
		\includegraphics[width=\linewidth]{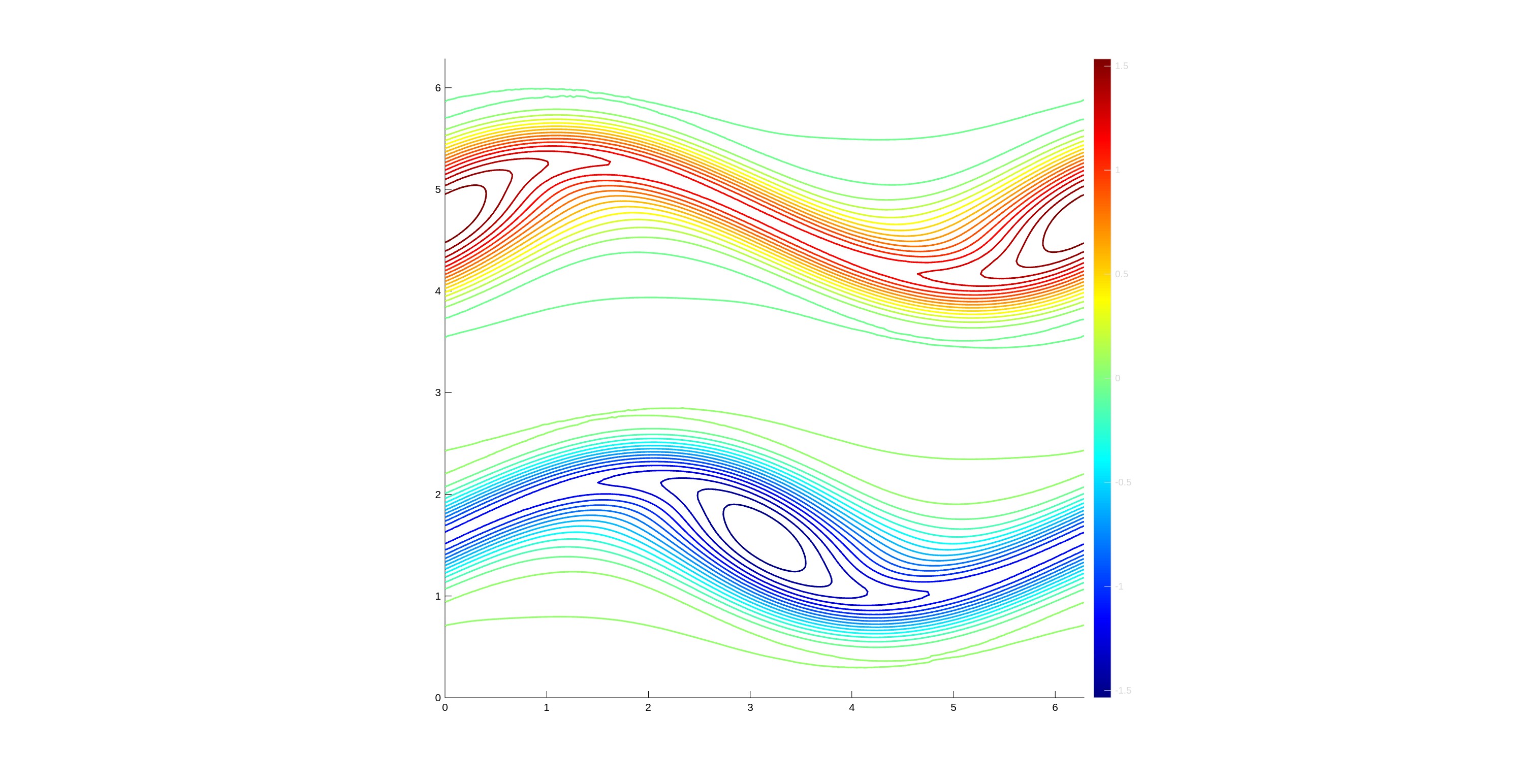}
		\caption{Semi-Lagrangian method, time 4.93993.}
		\label{fig:SLVorticity5}
	\end{subfigure}
	\begin{subfigure}[b]{0.49\textwidth}
		\centering
		\includegraphics[width=\linewidth]{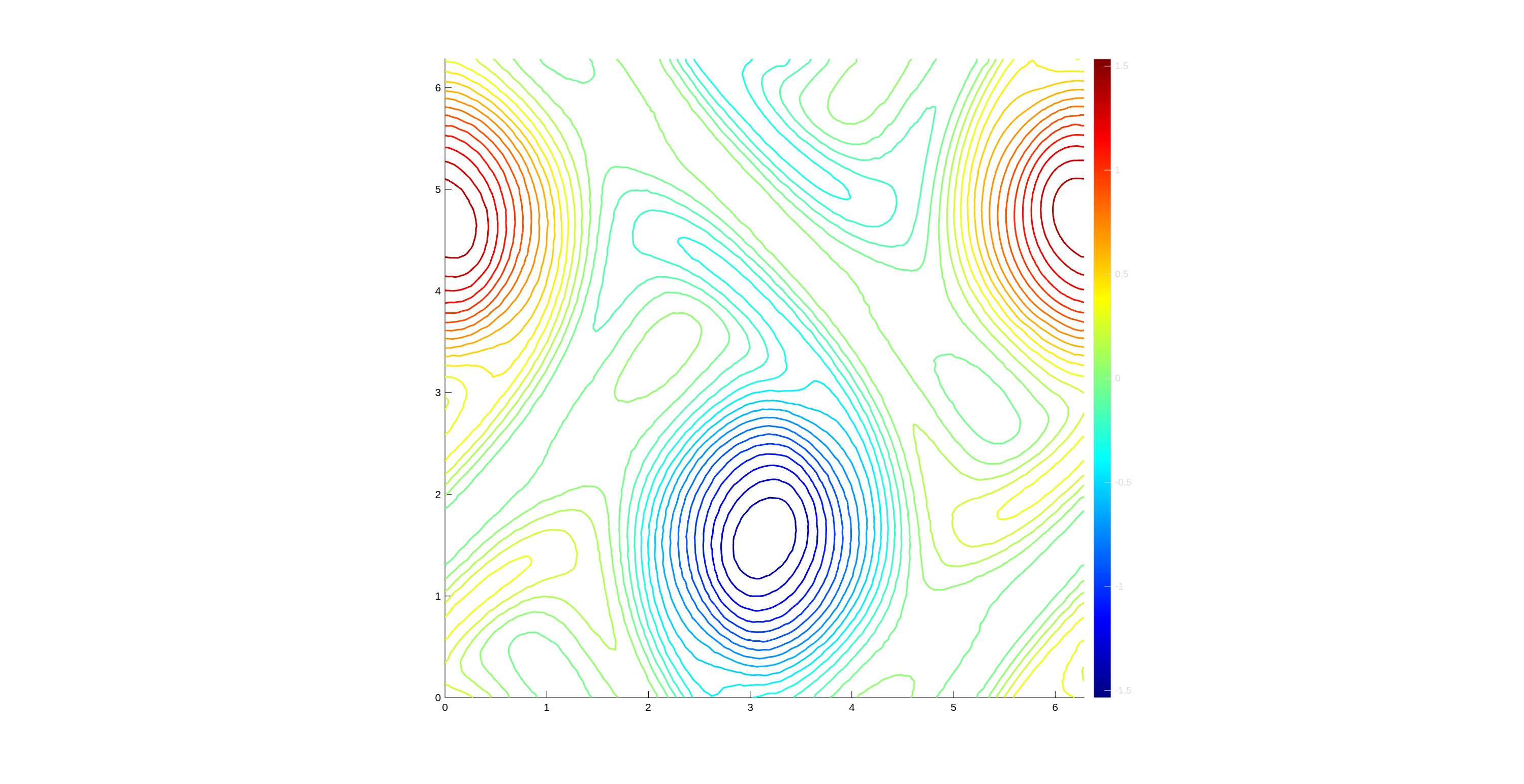}
		\caption{Meshfree MUSCL method, time 10.}
		\label{fig:MUSCLVorticity10}
	\end{subfigure}
	\begin{subfigure}[b]{0.49\textwidth}
		\centering
		\includegraphics[width=\linewidth]{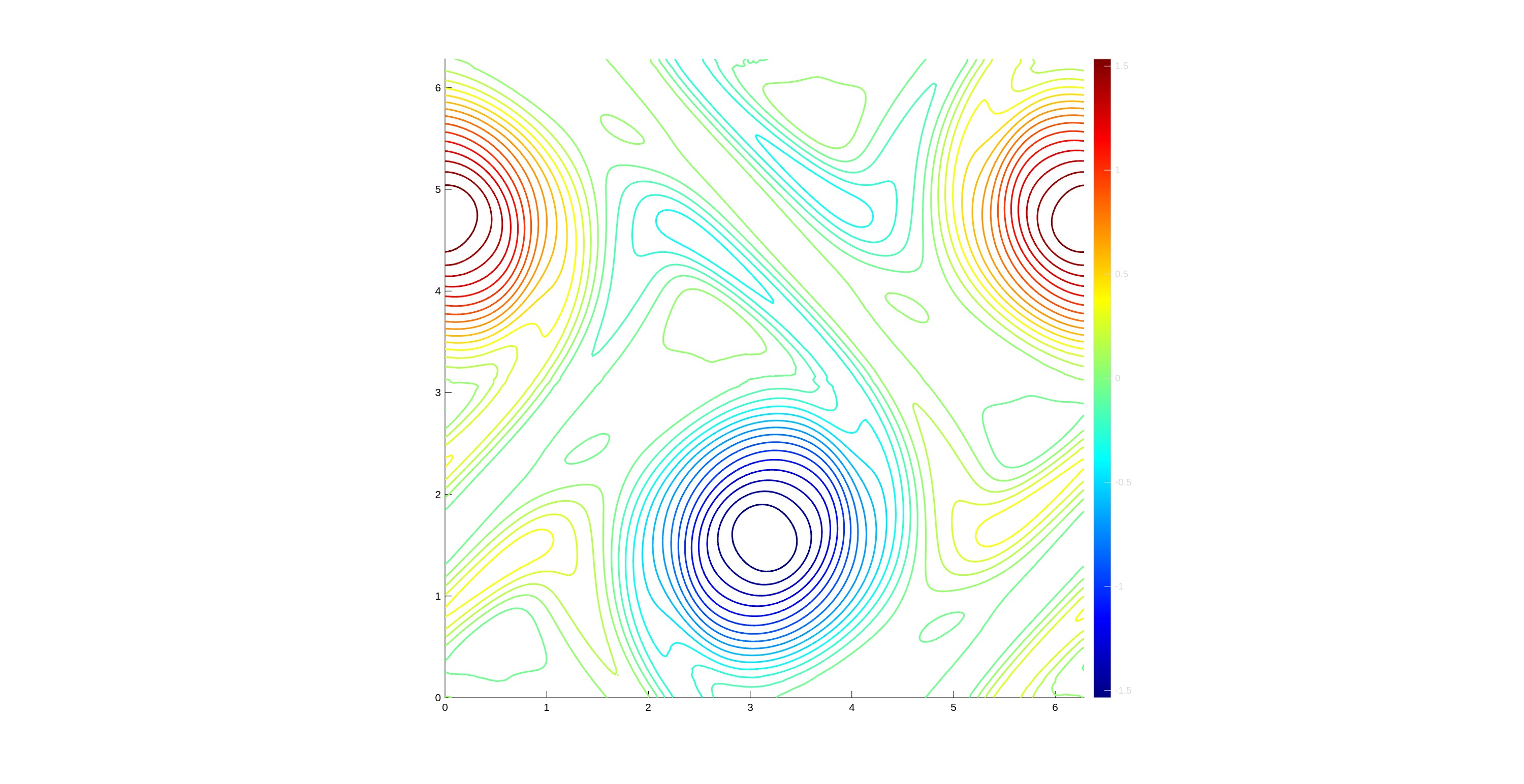}
		\caption{Semi-Lagrangian method, time 10.}
		\label{fig:SLVorticity10}
	\end{subfigure}
	\caption{Plots of the vorticity for the shear flow problem.}
	\label{fig:ShearFlowVorticity}
\end{figure}

The vorticity at time 4.93993 and 10 for both methods are plotted in figure \ref{fig:ShearFlowVorticity}. Using both algorithms, the flow clearly develops a vortex. We observe that the second-order in space meshfree method is only slightly more diffusive than the third-order in space semi-Lagrangian method. We conclude that the meshfree method is capable of resolving challenging flows. 

\section{Conclusion}
\label{section:conclusion}
We have presented a numerical method for simulating rarefied gases that interact with moving boundaries and rigid bodies. The gas is described by the BGK equation and solved using an Arbitrary Lagrangian-Eulerian method in which grid points move with the local mean velocity of the gas. This formulation allows for easy treatment of irregular domains and moving boundaries. We use higher-order meshless MUSCL schemes with MOOD to discretised the equation on the irregular, moving grid. Diffusive-reflective boundary conditions on the walls are treated using a new procedure that avoid an iterative solver or extrapolation procedure. We perform several tests and comparisons with different methods to validate the performance of the numerical method. In future work, we would like to tackle three-dimensional problems and mixtures of gases. 

\textbf{Acknowledgments: } This work was supported by the European Union’s Framework Program for Research and Innovation Horizon Europe under the Marie~Skłodowska-Curie Doctoral Networks action (HORIZON-MSCA-2021-DN-01), Grant Agreement No.~101072546 (DATAHYKING), and by the German Research Foundation (DFG) under Grant No.~KL~1105/35. G.~Russo additionally acknowledges the support of the Italian Ministry of University and Research (MUR) through the PRIN~2022 Project No.~2022KA3JBA, entitled ``Advanced numerical methods for time-dependent parametric partial differential equations with applications.''

\appendix

\section{The u2 detection criterion in two space dimensions}
\label{appendix:MOOD2D}
In two space dimensions, criteria \eqref{eq:MOODOscillationCheck} and \eqref{eq:MOODSignCheck} are changed to include checks on the oscillations indicators in the \(y\)-direction by \begin{align}
    \left( \mathcal{Y}_{i, k}^{\rm min} \mathcal{Y}_{i, k}^{\rm max} > - \delta \right) &\text{ and } \left( \mathcal{X}_{i, k}^{\rm min} \mathcal{X}_{i, k}^{\rm max} > - \delta \right) \\
    \left( \frac{\tilde{\mathcal{Y}}_{i, k}^{\rm min}}{\tilde{\mathcal{Y}}_{i, k}^{\rm max}} \geq 1/2 \text{ or } \tilde{\mathcal{Y}}_{i, k}^{\rm max} < \delta \right) &\text{ and } \left( \frac{\tilde{\mathcal{X}}_{i, k}^{\rm min}}{\tilde{\mathcal{X}}_{i, k}^{\rm max}} \geq 1/2 \text{ or } \tilde{\mathcal{X}}_{i, k}^{\rm max} < \delta \right). 
\end{align} The oscillation indicators in the \(y\)-direction are straightforward extensions of \eqref{eq:oscillationX1} and \eqref{eq:oscillationX2}, \begin{align}
\Tilde{\mathcal{Y}}_{i, k}^{\rm min} &= \min_{j \in \mathcal{C}_i} \left( \left| \pdv[2]{\Tilde{g}_{i, k}}{y} \right| , \left| \pdv[2]{\Tilde{g}_{j, k}}{y} \right| \right) \text{ and } \Tilde{\mathcal{Y}}_{i, k}^{\rm max} = \max_{j \in \mathcal{C}_i} \left( \left| \pdv[2]{\Tilde{g}_{i, k}}{y} \right|, \left| \pdv[2]{\Tilde{g}_{j, k}}{y} \right| \right), \\
\mathcal{Y}_{i, k}^{\rm min} &= \min_{j \in \mathcal{C}_i} \left( \pdv[2]{\Tilde{g}_{i, k}}{y}, \pdv[2]{\Tilde{g}_{j, k}}{y} \right) \text{ and } \; \mathcal{Y}_{i, k}^{\rm max} = \max_{j \in \mathcal{C}_i} \left( \pdv[2]{\Tilde{g}_{i, k}}{y}, \pdv[2]{\Tilde{g}_{j, k}}{y} \right). \end{align}
\newpage
\printbibliography

\end{document}